\documentclass[12pt,oneside]{amsart}
\pdfoutput=1 
\usepackage{geometry}
\usepackage{graphicx}
\usepackage{spectralsequences}
\usepackage{todonotes}
\usepackage{amsmath}
\usepackage{amsfonts}
\usepackage{amssymb}
\usepackage{tikz-cd}
\usepackage{bbm}
\usepackage{comment}
\usepackage{stmaryrd}
\usepackage{fullpage}
\usepackage{epstopdf}
\usepackage{scalerel}
\usepackage{amsthm}
\usepackage{mathrsfs} 
\usepackage{tikz}
\usepackage{dsfont}
\usetikzlibrary{matrix}
\usepackage{mathtools}
\usepackage[hidelinks]{hyperref}
\usepackage{cleveref}
\usepackage{xcolor}
\usepackage{enumitem}
\usepackage{tensor}
\usepackage[utf8]{inputenc}
\usepackage{csquotes}
\usepackage[english]{babel}
\usepackage{mymacros}

\hypersetup{
	colorlinks,
	linkcolor={red!50!black},
	citecolor={blue!50!black},
	urlcolor={blue!80!black}
}

\tikzset{
	rot90/.style={anchor=south, rotate=90, inner sep=.5mm}
}
\tikzset{
	rot45/.style={anchor=south, rotate=-45, inner sep=.5mm}
}

\setlist[itemize]{leftmargin=20pt}

\theoremstyle{definition}

\newtheorem{nul}{}[section]
\newtheorem{dfn}[nul]{Definition}

\newtheorem{rmk}[nul]{Remark}

\newtheorem{cnstr}[nul]{Construction}
\newtheorem{cnv}[nul]{Convention}
\newtheorem{ntn}[nul]{Notation}
\newtheorem{exm}[nul]{Example}

\newtheorem{rec}[nul]{Recollection}

\newtheorem{wrn}[nul]{Warning}

\newtheorem*{dfn*}{Definition}
\newtheorem*{axm*}{Axiom}
\newtheorem*{ntn*}{Notation}
\newtheorem*{exm*}{Example}
\newtheorem*{exr*}{Exercise}
\newtheorem*{int*}{Intuition}
\newtheorem*{qst*}{Question}
\newtheorem*{rmk*}{Remark}

\theoremstyle{plain}

\newtheorem{thm}[nul]{Theorem}
\newtheorem{prop}[nul]{Proposition}

\newtheorem{lem}[nul]{Lemma}

\newtheorem{cor}[nul]{Corollary}

\newtheorem*{thm*}{Theorem}
\newtheorem*{prop*}{Proposition}
\newtheorem*{cor*}{Corollary}
\newtheorem*{lem*}{Lemma}
\newtheorem*{cnj*}{Conjecture}

\AtBeginEnvironment{rmk}{\pushQED{\qed}}
\AtEndEnvironment{rmk}{\popQED\endrmk}
\AtBeginEnvironment{dfn}{\pushQED{\qed}}
\AtEndEnvironment{dfn}{\popQED\enddfn}
\AtBeginEnvironment{cnstr}{\pushQED{\qed}}
\AtEndEnvironment{cnstr}{\popQED\endcnstr}
\AtBeginEnvironment{rec}{\pushQED{\qed}}
\AtEndEnvironment{rec}{\popQED\endrec}
\AtBeginEnvironment{cnv}{\pushQED{\qed}}
\AtEndEnvironment{cnv}{\popQED\endcnv}
\AtBeginEnvironment{ntn}{\pushQED{\qed}}
\AtEndEnvironment{ntn}{\popQED\endntn}
\AtBeginEnvironment{exm}{\pushQED{\qed}}
\AtEndEnvironment{exm}{\popQED\endexm}
\AtBeginEnvironment{wrn}{\pushQED{\qed}}
\AtEndEnvironment{wrn}{\popQED\endwrn}

\title{Topological Hochschild homology of the image of j}
\newtoggle{draft}
\togglefalse{draft}
\makeatletter
\def\@tocline#1#2#3#4#5#6#7{\relax
	\ifnum #1>\c@tocdepth 
	\else
	\par \addpenalty\@secpenalty\addvspace{#2}%
	\begingroup \hyphenpenalty\@M
	\@ifempty{#4}{%
		\@tempdima\csname r@tocindent\number#1\endcsname\relax
	}{%
		\@tempdima#4\relax
	}%
	\parindent\z@ \leftskip#3\relax \advance\leftskip\@tempdima\relax
	\rightskip\@pnumwidth plus4em \parfillskip-\@pnumwidth
	#5\leavevmode\hskip-\@tempdima
	\ifcase #1
	\or\or \hskip 1em \or \hskip 2em \else \hskip 3em \fi%
	#6\nobreak\relax
	\hfill\hbox to\@pnumwidth{\@tocpagenum{#7}}\par
	\nobreak
	\endgroup
	\fi}
\makeatother

\iftoggle{draft} {
							\usepackage[notcite]{showkeys}
                         	\newcommand{\NB}[1]{\todo[color=gray!40]{#1}}
                         	\newcommand{\TODO}[1]{\todo[color=red]{#1}}
}{ 
                         	\newcommand{\NB}[1]{}
                         	\newcommand{\TODO}[1]{}
                         	\renewcommand{\todo}[1]{}
                         	\renewcommand{\todo}[1]{}
}
\date{\today}

\author{David Jongwon Lee}
\address{Department of Mathematics, Northwestern University, Evanston, IL, USA}
\email{davidlee@northwestern.edu}

\author{Ishan Levy}
\address{Department of Mathematics, Institute for Advanced Study, Princeton, NJ, USA}
\email{ishanl@ias.edu}

\begin{document}
\begin{abstract}
	We compute the mod $(p,v_1)$ and mod $(2,\eta,v_1)$ $\THH$ of many variants of the image-of-$J$ spectrum. In particular, we do this for $j_{\zeta}$, whose $\TC$ is closely related to the $K$-theory of the $K(1)$-local sphere. We find in particular that the failure for $\THH$ to satisfy $\ZZ_p$-Galois descent for the extension $j_{\zeta} \to 
	\ell_p$ corresponds to the failure of the $p$-adic circle to be its own free loop space. For $p>2$, we also prove the Segal conjecture for $j_{\zeta}$, and we compute the $K$-theory of the $K(1)$-local sphere in degrees $\leq 4p-6$.
\end{abstract}
\maketitle
\setcounter{tocdepth}{2}
\tableofcontents
\section{Introduction}

The algebraic $K$-theory of the $K(1)$-local sphere, or $K(L_{K(1)}\SP)$, is an object capturing fundamental structural information about the $K(1)$-local category. Part of Ausoni--Rognes' original vision of chromatic redshift was that it could be understood, at least $T(2)$-locally, via Galois hyperdescent. More specifically, they conjectured \cite[pg 4]{ausonirognes} that the map
\[
K(L_{K(1)}\SP)\otimes V \to K(\KU_p)^{h\ZZ_p^\times}\otimes V
\]
is an equivalence in large degrees when $V$ is a type $2$ finite spectrum. The $T(n+1)$-local $K$-theory of Morava $E$-theory has been shown in \cite[Theorem 1.10]{clausen2016descent} to have Galois descent for finite subgroups of the Morava stabilizer group. Moreover, recent work of Ben Moshe--Carmeli--Schlank--Yanovski \cite{cyclotomicredshift} shows that $L_{K(2)}K(L_{K(1)}\SP) \to L_{K(2)}K(\KU_p)^{h\ZZ_p^{\times}}$ is an equivalence, i.e that Galois hyperdescent is satisfied $K(2)$-locally.

Recent work of the second author \cite{levy2022algebraic} has made $K(L_{K(1)}\SP)$ an integrally accessible object. If we consider the connective Adams summand $\ell_p$ (or $\ko_2$ for $p=2$) as a $\ZZ$-equivariant $\EE_{\infty}$-ring via the Adams operation $\Psi^{1+p}$, then $j_{\zeta}$ is defined to be its $\ZZ$-homotopy fixed points. Then it is shown that there is a cofiber sequence
\[
K(j_{\zeta}) \to K(L_{K(1)}\SP) \to \Sigma K(\FF_p)
\]
split on $\pi_*$. It is also shown that the Dundas--Goodwillie--McCarthy square

\begin{center}
	\begin{tikzcd}
		K(j_{\zeta}) \ar[r]\ar[d] & \TC(j_{\zeta})\ar[d]\\
		K(\ZZ_p)\ar[r] & \TC(\ZZ_p^{B\ZZ})
	\end{tikzcd}
\end{center}
is a pullback square. 
The three spectra $K(\FF_p), K(\ZZ_p)$, and $\TC(\ZZ_p^{B\ZZ})$\footnote{$\ZZ_p^{B\ZZ}$ denotes the cochains on the circle $B\ZZ$ with coefficients in $\ZZ_p$. Its $\TC$ is essentially the nil-$\TC$ of $\ZZ_p$ by \cite[Corollary 4.5]{land2023k}, which is studied in \cite{hesselholt2004rham}.} are understood, so understanding $K(L_{K(1)}\SP)$ is essentially reduced to understanding $\TC(j_{\zeta})$.

The primary goal of this paper is to understand $\THH(j_{\zeta})$ modulo $(p,v_1)$ and $(2,\eta,v_1)$, which is the first step in understanding $\TC(j_{\zeta})$.

\begin{thm}\label{thm:main}
	For $p>2$, there is an isomorphism of rings
    \[
    \pi_*\THH(j_{\zeta})/(p,v_1) \cong 
	\pi_*\THH(\ell_p)/(p,v_1)\otimes_{\mathbb F_p} \mathrm{HH}_*(\FF_p^{B\ZZ}/\FF_p).
    \]
	For $p=2$, there is an isomorphism of rings 
    \[
    \pi_*\THH(j_{\zeta})/(2,\eta,v_1) \cong 
	\pi_*\THH(\ko_2)/(2,\eta,v_1)\otimes_{\mathbb F_2} \mathrm{HH}_*(\FF_2^{B\ZZ}/\FF_2).
    \]
\end{thm}

Each of the terms on the right hand side of the equivalences is well understood. The ring $\pi_*\THH(\ell_p)/(p,v_1)$ can be found in \cite{mcclure1993topological} or \Cref{exm:THHell}, and $\pi_*\THH(\ko_2)/(2,\eta,v_1)$ can be found in \cite{angeltveit2005hopf} or \Cref{exm:thhko}.
The last tensor factor is given in \Cref{lem:spheretriv} as
\[
\HoH_*(\FF_p^{B\ZZ}/\FF_p) \cong \Lambda[\zeta]\otimes \ctf{\ZZ_p}
\]
where $|\zeta| = -1$ and $\ctf{\ZZ_p}$ denotes the ring of continuous functions from $\ZZ_p$ to $\FF_p$.

The factor of $\ctf{\ZZ_p}$ can be viewed as the failure of descent at the level of $\THH$ for the $\ZZ_p$-Galois extension coming from the $\ZZ_p$-action on $\ell_p$ and $\ko_2$. More precisely, at the level of $\pi_*$, the natural map
\[
\THH(\ell_p^{h\ZZ})/(p,v_1) \to \THH(\ell_p)^{h\ZZ}/(p,v_1)
\]
is induced by the base change along $\ctf{\ZZ_p} \to \FF_p$ that sends a continuous function to its value at $0$ (\Cref{rmk:THHellhZ}).

This phenomenon can be explained by interpreting $\THH$ in terms of free loop spaces. If $X$ is a pro-$p$-finite space, then the $\FF_p$-Hochschild homology of the cochain algebra $C^*(X;\FF_p)$ is computed as
\[
	\HoH(C^\ast(X;\mathbb F_p)/\mathbb F_p) = C^\ast(LX;\mathbb F_p)
\]
where $LX$ is the free loop space of $X$. Since $C^*(B\ZZ_p;\FF_p) \cong \FF_p^{B\ZZ}$, the failure of the descent
\[
	\HoH(\mathbb F_p^{B\ZZ}/\mathbb F_p) \not\simeq \HoH(\mathbb F_p/\mathbb F_p)^{B\ZZ}
\]
is explained by the fact that $B\ZZ_p$ is not $LB\ZZ_p \cong B\ZZ_p\times \ZZ_p$. For any $p$-complete $\EE_{\infty}$-ring $R$ with a trivial $\mathbb Z$-action, this completely accounts for the failure of $p$-complete $\THH$ to commute with $\ZZ$-fixed points (\Cref{cor:trivfixed}). The content of \Cref{thm:main} is that the same phenomenon happens for $\THH(j_{\zeta})$ on $\pi_*$ mod $(p,v_1)$ or $(2,\eta,v_1)$, even though the action is no longer trivial. In particular, \Cref{thm:main} implies that there is an isomorphism of rings
\[
\pi_*\THH(j_{\zeta})/(p,v_1) \cong \pi_*\THH(\ell_p^{B\ZZ})/(p,v_1)
\]
where $\ell_p^{B\ZZ}$ is the homotopy fixed points of $\ell_p$ by a \textit{trivial} $\ZZ$-action.

The key idea in our proof of \Cref{thm:main} is to run the spectral sequence for $\THH$ obtained by filtering $j_{\zeta}$ via the homotopy fixed point filtration, and showing that the differentials in the associated spectral sequence behave similarly enough to the case of a trivial action. To understand the associated graded algebra of the homotopy fixed point filtration, we further filter it by the $p$-adic filtration. At the level of the associated graded of both filtrations, $j_{\zeta}$ is indistinguishable from the fixed points by a trivial action, and we show that mod $(p,v_1)$ and $(2,\eta,v_1)$ this remains true at the level of homotopy rings after running the spectral sequences for $\THH$ of those filtrations.

The phenomenon that the $\ZZ$-action on $\ell_p$ behaves like the trivial one is shown in \cite{telescope} to asymptotically hold even at the level of cyclotomic spectra. More precisely, it is shown there that given any fixed type $3$ finite spectrum $V$, for all sufficiently large $k$,
\[
\THH(\ell_p^{hp^k\ZZ})\otimes V \cong \THH(\ell_p^{B\ZZ})\otimes V
\]
as cyclotomic spectra.

It is shown then that the failure of descent we observe on $\THH$ continues at the level of the $T(2)$-local $\TC$. Combining this with the aforementioned hyperdescent result of the $K(2)$-local $K$-theory and the formula for the $K$-theory of the $K(1)$-local sphere, this implies that $L_{T(2)}K(L_{K(1)}\mathbb S)$ is not $K(2)$-local and hence is a counterexample to the height $2$ telescope conjecture. 
In particular, this implies that the map
\[
K(L_{K(1)}\SP)\otimes V \to K(\KU)^{h\ZZ_p^\times}\otimes V
\]
considered by Ausoni--Rognes is \textit{not} an equivalence in large degrees.

The ring $\ctf{\ZZ_p}$ that appears in our formula for $\THH(j_{\zeta})$ is a key ingredient in \cite{telescope} to maintain asymptotic control over $\THH(j_{\zeta,k})$ as a cyclotomic spectrum, and is one of the advantages of $j_{\zeta}$ versus the usual connective image-of-$J$ spectrum $j = \tau_{\geq0}j_{\zeta}$. If one was only interested in understanding $L_{T(2)}K(L_{K(1)}\SP)$, there are isomorphisms
\[
L_{T(2)}K(L_{K(1)}\SP) \cong L_{T(2)}\TC(j) \cong L_{T(2)}\TC(j_{\zeta})
\]
so one can in principle approach the telescopic homotopy via $\TC(j)$ instead of $\TC(j_{\zeta})$.

However, $j$ is not as well behaved as $j_{\zeta}$ is, as we now explain. We extend our methods for computing $\THH(j_{\zeta})$ in \Cref{sec:j,sec:finiteextn} to compute $\THH$ of $j$, giving a relatively simple proof of the result below due to Angelini-Knoll and H\"oning \cite{angelini2021topological,honing2021topological}.

\begin{thm}\label{thm:j}
	For $p>3$, the ring\footnote{We also compute an associated graded ring $\THH(j)/(p,v_1)$ for $p=3$ (see \Cref{thm:thhj}), but are unable to solve multiplicative extension problems coming from the fact that $j/(p,v_1)$ is not an associative algebra for $p=3$. Nonassociative multiplicative extensions aren't considered in \cite{angelini2021topological}, so the results of that paper also only compute an associated graded ring for $p=3$.} $\pi_*\THH(j)/(p,v_1)$ is the homology of the CDGA
	\[
	\mathbb F_p[\mu_2]\otimes\Lambda[\alpha_1,\lambda_2,a]\otimes\Gamma[b],\quad d(\lambda_2)=a\alpha_1
	\]
	\[ |b| = 2p^2-2p , \; |a| = 2p^2-2p-1,\; |\lambda_2| = 2p^2-1,\; |\mu_2| = 2p^2\]
	For $k\geq 1$ and any $p>2$, we have an isomorphism of rings $$\pi_*\THH(\tau_{\geq0}(\ell_p^{hp^k\ZZ}))/(p,v_1) \cong \pi_*\THH(\ell_p)/(p,v_1)\otimes \HoH_*(\tau_{\geq0}\FF_p[v_1]^{Bp^k\ZZ}/\FF_p[v_1])/v_1.$$
\end{thm}

The ring $\HoH_*(\tau_{\geq0}\FF_p[v_1]^{Bp^k\ZZ}/\FF_p[v_1])/v_1$ is described in \Cref{prop:hhjgr/p}: it is isomorphic\footnote{By rescaling, there is an isomorphism $p^k\ZZ \cong \ZZ$, so this ring doesn't depend on $k$.} to $\Gamma[d\alpha_{1/p^k}]\otimes \Lambda_{\FF_p}[\alpha_{1/p^k}]$ where $\alpha_{1/p^k}$ is a class in degree $2p-3$ and $d\alpha_{1/p^k}$ is a divided power generator in degree $2p-2$.

In the above theorem, $\pi_*\THH(j)/(p,v_1)$ is \textit{not} what one would expect in the case of the trivial action: there are two more differentials in the spectral sequence for the filtration we use to prove \Cref{thm:j} than what one would find for the trivial action. The differentials witness the fact that $\lambda_1,\lambda_2 \in \pi_*\THH(\ell_p)/(p,v_1)$ (see \Cref{exm:THHell}) don't lift to $\THH(j)/(p,v_1)$. Whereas most computations of $\THH$ in this paper use B\"okstedt's computation of $\THH(\FF_p)$ as their fundamental input, these differentials ultimately come from the Adams--Novikov spectral sequence.

A key difference between the $\THH$ of $j_{\zeta}$ and $j$ is that the ring $\ctf{\ZZ_p}$ that appeared in $\pi_*\THH(j_{\zeta})/(p,v_1)$ is replaced by a divided power algebra for $j$. The advantage of the ring $C^0(\ZZ_p;\FF_p)$ over a divided power algebra is that up to units, it consists entirely of idempotents, which decompose $\THH(j_{\zeta})$ as an $S^1$-equivariant spectrum into a continuous $\ZZ_p$-indexed family of spectra. This decomposition is not evidently present in $\THH(j)$.

Another advantage of $j_{\zeta}$ over $j$ is that $j_{\zeta}$ satisfies the $\THH$ Segal conjecture while $j$ doesn't, which we show for $p>2$ in \Cref{sec:seg}:

\begin{thm}\label{thm:segintro}
	For $p>2$, the cyclotomic Frobenius map
	\[
	\THH(j_{\zeta})/(p,v_1) \to \THH(j_{\zeta})^{tC_p}/(p,v_1)
	\]
	has $(2p-3)$-coconnective fiber, but the fiber of the cyclotomic Frobenius map
	\[
	\THH(j)/(p,v_1) \to \THH(j)^{tC_p}/(p,v_1)
	\]
	is not bounded above.
\end{thm}

The Segal conjecture for a ring $j$ is a necessary condition \cite[Proposition 2.25]{antieau2021cartier} for the Lichtenbaum--Quillen conjecture to hold, i.e for $\TR(j)\otimes V$ to be bounded above for any finite type $3$ spectrum $V$. Thus \Cref{thm:segintro} implies that $j$ doesn't satisfy the Lichtenbaum--Quillen conjecture. On the other hand, \Cref{thm:segintro} is a key ingredient in proving the Lichtenbaum--Quillen conjecture for $j_{\zeta}$ as carried out in \cite{telescope}. This Lichtenbaum--Quillen conjecture can be viewed as the part of Ausoni--Rognes's conjecture that is true. Namely, it implies that the map 
$$K(L_{K(1)}\SP)\otimes V \to K(L_{K(1)}\SP)\otimes V[v_2^{-1}]$$
is an equivalence in large degrees for $V$ a type $2$-complex. The Lichtenbaum--Quillen conjecture for $\ell^{hp^k\ZZ}$ for $k\gg0$ is a key ingredient in \cite{telescope} to show that the telescope conjecture fails at height $2$.

In \Cref{sec:stable}, we show how $\THH$ computations can give information about $\TC$ in the stable range. For a map of $\EE_1$-rings $f:R \to S$, recall that the $\EE_1$-cotangent complex $L_{S/R}$ is the fiber of the multiplication map $S\otimes_RS \to S$ as an $S$-bimodule. We prove the following result:

\begin{thm}\label{thm:stablerange}
	Given a map of $\EE_1$-ring spectra $f:R \to S$, there is a natural map
	\[
	\fib\TC(f) \to \THH(S;L_{S/R}).
	\]If $f$ is an $n$-connective map of $(-1)$-connective rings for $n\geq 1$, this natural map is $(2n+1)$-connective.
\end{thm}

A consequence of \Cref{thm:stablerange} is that the natural map above can be identified with the linearization map in the sense of Goodwillie calculus for the functor $\fib\TC(f): \Alg(\Sp)_{/S} \to \Sp$ when $S$ is $(-1)$-connective and $f:R \to S$ is $1$-connective.

In the case the map $f$ is a trivial square zero extension of connective rings, a $K$-theory version of the result was obtained as \cite[Theorem 3.4]{dundas1994stable}, and a $\TC$ version is essentially \cite[Theorem 4.10.1]{raskin2018dundas}\footnote{See also \cite{hesselholt1994stable} and \cite{lindenstrauss2012taylor}.}. The point of \Cref{thm:stablerange} is to have a version of the result that works for arbitrary maps of $\EE_1$-rings rather than trivial square-zero extensions, and for $(-1)$-connective rings instead of connective rings.

We use \Cref{thm:stablerange} to reprove basic facts about $\TC$, such as the understanding of the map $\TC(\SP_p) \to \TC(\ZZ_p)$ on $\pi_{2p-1}$. This is an ingredient in the computation of $\TC(\ZZ_p)$ as a spectrum (see \cite[Section 9]{bokstedt1993topological}).

We also apply \Cref{thm:stablerange} to compute the fiber of the map $\TC(j_{\zeta}) \to \TC(\ZZ_p^{B\ZZ})$ in the stable range, giving information about $K(L_{K(1)}\SP)$:

\begin{thm}\label{thm:kk1localstableintro}
	For $p>2$, there are isomorphisms $$\tau_{\leq 4p-6}\fib(\TC(j_{\zeta}) \to \TC(\ZZ_p^{B\ZZ})) \cong \Sigma^{2p-2}\ctf{\ZZ_p}$$ and $$K_{*}L_{K(1)}\SP \cong K_{*-1}\FF_p \oplus K_{*}\SP_p \oplus \pi_*\Sigma^{2p-2}\ctf{\ZZ_p}/\FF_p,\;\; *\leq 4p-6.$$
\end{thm}

In particular, for $p>2$, the infinite family of classes in the fiber of $\TC(j_{\zeta}) \to \TC(\ZZ_p^{B\ZZ})$ found in \cite{levy2022algebraic} are simple $p$-torsion, and completely account for all the classes in the stable range.

\subsection*{Acknowledgements}
We are very grateful to Robert Burklund, Sanath Devalapurkar, Jeremy Hahn, Mike Hopkins, Tomer Schlank, and Andy Senger for conversations related to this work. The second author is supported by the NSF Graduate Research Fellowship under Grant No. 1745302.

\subsection*{Notations and conventions}
\begin{itemize}
	\item The term category will refer to an $\infty$-category as developed by Joyal and Lurie.
	\item We refer the reader to \cite{nikolaus2018topological} for basic facts about $\THH$, which we freely use.
	\item For $\cC$ a monoidal category, and $R$ an $\EE_1$-algebra, we use $\THH_{\cC}(R)$ to denote the $\THH$ of $R$ in $\cC$. For $\cC = \Mod(S)$ for an $\EE_2$-algebra $S$, we also denote this by $\THH(R/S)$.
	\item $\Map(a, b)$ will denote the space of maps from 
	$a$ to $b$ (in some ambient category).
	\item Tensor products and $\THH$ are implicitly $p$-completed.
	\item We use $\Lambda[x]$ and $\Gamma[x]$ to denote exterior and divided power algebras in homotopy rings.
	\item In an $\FF_p$-vector space, we use $a \doteq b$ to mean that $a = cb$ for some unit $c \in \FF_p^\times$, and $a\dotmapsto b$ to mean that $a$ is sent to $b$ up to a unit in $\FF_p^\times$.
	\item Conventions about filtrations and spectral sequences are addressed in \Cref{sec:filtration}.
	\item For a pro-finite set $A$, we use $C^0(A;\FF_p)$ to denote the ring of continuous functions from $A$ to $\FF_p$.
	\item Let $\mathcal{D}$ be a monoidal category acting on a category $\mathcal{C}$. Given objects $X \in \mathcal{C}, Z \in \mathcal{D}$ with a self map $f:X\otimes Z\to X$, we use $X/f$ to denote the cofibre of this map. We use $X/(f_1,\dots,f_n)$ to denote $(\dots(X/f_1)/\dots)/f_n$, where each $f_i$ is a self map of $X/(f_1,\dots,f_{i-1})$.
\end{itemize}

\section{Filtrations}\label{sec:filtration}
In this section, we set up notation for working with filtered objects and explain how to put filtrations on $\ell_p$, $\ko_2$, $j_{\zeta}$, and $j$, as well as for finite extensions. Our constructions amount to the filtration coming from the homotopy fixed point spectral sequences computing those objects, which in all cases except for $j_{\zeta}$, is also the Adams--Novikov filtration.

\subsection{Filtered objects and spectral sequences}
Let $\mathcal C$ be a presentably symmetric monoidal stable category with accessible  
$t$-structure compatible with the symmetric monoidal structure. Let $\Fil(\mathcal C) = \Fun(\ZZ_{\leq}^{\mathrm{op}},\mathcal C)$ be the category of decreasingly filtered objects, and let $\gr(\mathcal C) = \Fun(\ZZ,\mathcal C)$ be the category of graded objects, so that both are symmetric monoidal via Day convolution. Basic properties of these categories are developed in \cite{rotation} and \cite[Appendix B]{burklund2020}. Given an object $x \in \Fil(\mathcal C)$ or $\gr(\mathcal C)$, we write $x_i$ for the value at $i \in \ZZ$. The left adjoint of the functor $(-)_i$ in the case of $\Fil(\mathcal C)$ is the functor $(-)^{0,i}$, defined for $c\in \mathcal C$ by
\[
(c^{0,i})_j = \begin{cases}
c&(j\leq i)\\
0&(j>i).
\end{cases}
\]We also use the notation (where $\Sigma$ denotes suspension in $\mathcal{C}$)\todo{this notation isn't used consistently throughout the paper.}
\begin{align*}
c^{k,n} &:= \Sigma^kc^{0,n+k}\\
\pi_{k,n}^{\heart}x &:= \pi_{k}^{\heart}x_{n+k}\\
\pi_{k,n}x&:= \pi_kx_{n+k} = \pi_0\Map(\unit^{k,n},x)
\end{align*}
and use $c$ to also denote $c^{0,0}$. There is class $\tau\in \pi_{0,-1}\unit^{0,0}$ called the \emph{filtration parameter} which is adjoint to the identity map $\unit \to (\unit^{0,0})_{-1}$.\footnote{This usage of $\tau$ is compatible with the usage in \cite[Appendix B]{burklund2020}.}

The functor $(-)^{0,0}: \mathcal{C} \to \Fil(\mathcal{C})$ is a symmetric monoidal fully faithful functor, which we refer to as the \textit{trivial filtration}. We often identify an object $c \in \mathcal{C}$ with the trivial filtered object in $\Fil(\mathcal{C})$.
 In fact, $\gr(\mathcal C)$ can be identified with $\Mod_{\cof \tau}(\Fil(\mathcal C))$, so that taking associated graded amounts to base changing to $\cof \tau$. Given an object $x \in \Fil(\mathcal C)$, we let $\gr x\in\gr(\mathcal C)$ denote the associated graded object, so that $(\gr x)_i = \gr_ix = \cof(x_{i+1} \xrightarrow{\tau} x_i)$.

On the other hand, there is an identification $\Fil(\mathcal{C})[\tau^{-1}] \cong \mathcal C$, so that given a filtered object $x \in \Fil(\mathcal{C})$, its underlying object $ux \in \mathcal{C}$, given by $\colim_i x_i$, is identified with $x[\tau^{-1}]$. Under the assumption that the $t$-structure is compatible with filtered colimits, we have an isomorphism $\pi_{**}^{\heart}x[\tau^{-1}] \cong \pi_{*}^{\heart}ux\otimes \ZZ[\tau^{\pm1}]$.
\begin{cnstr}
Given a filtered object $x\in\Fil(\mathcal C)$, there is a spectral sequence which we refer to as \textit{the spectral sequence associated with $x$}.
\[
	E_1^{s,t} = \pi^{\heart}_{t-s,s}\gr x=\pi^{\heart}_{t-s}(\gr x)_t\implies \pi^{\heart}_{t-s}(ux)
\]
The $d_r$-differential is a map from $E_r^{s,t}$ to $E_r^{s+r+1,t+r}$, which is a page off from the usual Adams convention, i.e. our $d_r$ differential would be the $d_{r+1}$ differential in the Adams convention. We shall say \emph{Adams weight} and \emph{filtration degree} to refer to the bidegrees $s$ and $t$, respectively.
\end{cnstr}

In addition to the spectral sequence associated with $x$, there is also the $\tau$-Bockstein spectral sequence, which has signature
	$$E_1^{**} = (\pi_{**}^{\heart}\gr x)[\tau] \implies \pi_{**}^{\heart}x$$

We do not use the following lemma, but we state it as an exercise to help acquaint the unfamiliar reader with filtered objects. The $\tau$-inverted $\tau$-Bockstein spectral sequence refers to the spectral sequence obtained from the $\tau$-Bockstein spectral sequence by inverting $\tau$ on each page.
\begin{lem} Let $x \in \Fil(\mathcal C)$. For each $r\geq1$, the $E_r$-page of the $\tau$-inverted $\tau$-Bockstein spectral sequence for $x$ is isomorphic to $\ZZ[\tau^{\pm}]$ tensored with the $E_r$-page of the spectral sequence associated with $x$. Moreover, the $d_r$ differential on the former is given by $\tau^r$ times the $d_r$ differential on the latter. The filtration on $\pi^{\heart}_{**}x[\tau^{\pm1}]$ coming from the spectral sequence agrees with the filtration on $\pi^{\heart}_{*}x\otimes \ZZ[\tau^{\pm}]$ coming from the filtration on $x$.
\end{lem} 
\begin{proof}
	These statements can be checked for example by using explicit formulas for the pages and differentials. See, for example, \cite[Construction 1.2.2.6]{HA}.
\end{proof}

\subsection{$t$-structures}
We turn to studying $t$-structures on categories of filtered objects.
Our ability to produce $t$-structures comes from the following general result.
\begin{lem}[{\cite[Proposition 1.4.4.11]{HA}}]\label{lem:generalbigt}
	Let $\mathcal C$ be a presentable stable category.
	If $\{X_\alpha\}$ is a small collection of objects in $\mathcal C$, then there is an accessible $t$-structure $(\mathcal C_{\geq 0}, \mathcal C_{\leq 0})$ on $\mathcal C$ such that $\mathcal C_{\geq 0}$ is the smallest full subcategory of $\mathcal C$ containing each $X_\alpha$ and closed under colimits and extensions. The full subcategory of coconnective objects is characterized by the condition that $Y \in \mathcal C_{\leq 0}$ if and only if $\Map(\Sigma X_{\alpha},Y) = 0$ for each $X_\alpha$.
\end{lem}

\begin{dfn}\label{dfn:filtstr}
	Let $f:\ZZ \to \ZZ$ be a function. Define a $t$-structure $(\Fil(\mathcal C)^{f}_{\geq0},\Fil(\mathcal C)^{f}_{\leq0})$ on the underlying category $\Fil(\mathcal C)$ be the $t$-structure whose connective objects are generated by the objects $\Sigma^{f(i)}c^{0,i}$ for $c \in \mathcal C_{\geq0}$ and $i \in \ZZ$. We let $\tau^{f}_{\geq i}$ and $\tau^{f}_{\leq i}$ denote the associated truncation functors. We similarly define a $t$-structure $(\gr(\mathcal C)^f_{\geq0}, \gr(\mathcal C)^f_{\leq0})$ by taking the image of those objects under the functor $\gr$ to be the generators.
\end{dfn}

\begin{lem}\label{lem:ttruncate}
	Let $x \in \Fil(\mathcal C)$.
\begin{enumerate}
	\item $x \in \Fil(\mathcal C)^{f}_{\leq0}$ if and only if $x_i$ is $f(i)$-coconnective in $\mathcal C$ for each $i$.
	
	\item If $f$ is nondecreasing, then $x\in\Fil(\mathcal C)^{f}_{\geq0}$ iff $x_i$ is $f(i)$-connective for each $i$. In this case, the truncation functor $\tau^{f}_{\geq 0}$ is given by $(\tau^{f}_{\geq 0}x)_i = \tau_{\geq f(i)}(x_i)$.
	
	\item The same results hold for $(\gr(\mathcal C)^{f}_{\geq0},\gr(\mathcal C)^{f}_{\leq0})$.
\end{enumerate}
\end{lem}
\begin{proof}
	We prove the result for $\Fil(\mathcal C)$, as the result for $\gr(\mathcal C)$ is similar but easier. Coconnectivity can be checked by mapping in the generators of $\Fil(\mathcal C)^{f}_{\leq0}$. Because of the adjunction defining the functor $(-)^{0,n}$, the condition for coconnectivity follows.
	
	Now suppose $f$ is nondecreasing. To prove the claims, it suffices to show that if $x\in \Fil(\mathcal C)$ has $x_i \in \mathcal C_{\geq f(i)}$, then $x$ admits no maps to a coconnected object. If $y$ is a coconnected object, then $x_i$ admits no maps to $y_j$ for $j \leq i$ because $y_j$ is $f(j)$-coconnected, and since $f$ is nondecreasing, it is $f(i)$-coconnected. It follows that there are no nonzero maps of filtered objects $x \to y$.
\end{proof}

\begin{lem}\label{lem:tsymmon}
	The $t$-structures $\Fil(\mathcal C)^{f}, \gr(\mathcal C)^{f}$ are compatible with the symmetric monoidal structure if $f(0) = 0$ and $f(i) + f(j) \geq f(i+j)$.

\end{lem}
\begin{proof}
	The condition $f(0) = 0$ guarantees that the unit is connective.
	One needs to check that the tensor product of any pair of generators of $\Fil(\mathcal C)^{f}_{\geq0}$ is still in $\Fil(\mathcal C)^{f}_{\geq0}$. But the tensor product of $\Sigma^{f(i)}c^{0,i}$ and $\Sigma^{f(j)}d^{0,j}$ is $\Sigma^{f(i)+f(j)}(c\otimes d)^{0,i+j}$, which is in $\Fil(\mathcal C)^{f}_{\geq0}$ because $c\otimes d$ is in $\mathcal C_{\geq0}$ and so the assumption on $f$ shows that this is connective.
\end{proof}

The functor $\gr$ is right $t$-exact 
with respect to the $t$-structure corresponding to a nondecreasing function $f$, but not in general $t$-exact. In the following situation it preserves $\tau_{\geq0}$.

\begin{lem}\label{lem:grtruncate}
	Suppose that $c \in \Fil(\mathcal C)$, $f:\ZZ \to \ZZ$ is nondecreasing, $\pi^{\heart}_{k,i-k}c = 0$ for $f(i-1)\leq k < f(i)$, and $\pi^{\heart}_{f(i)-1,i-f(i)+2}c$ contains no simple $\tau$-torsion. Then $\tau_{\geq0}^f\gr(c) \cong \gr(\tau_{\geq0}^f(c))$ and $\tau_{\leq0}^{f}\gr(c) \cong \gr(\tau_{\leq0}^f(c))$.
\end{lem}

\begin{proof}
	It suffices to prove the statement for $\tau_{\geq0}$ since $\gr$ is exact. There is a cofiber sequence $c_{i+1} \xrightarrow{\tau} c_i \to \gr_ic$. By \Cref{lem:ttruncate} we would like $\tau_{\geq f(i+1)}c_{i+1} \to \tau_{\geq f(i)}c_i \to \tau_{\geq f(i)}\gr_ic$ to remain a cofiber sequence. From the exact sequence of homotopy groups, we see that we would like
	$\tau_{\geq f(i+1)}c_{i+1} = \tau_{\geq f(i)}c_{i+1}$ and $\pi^{\heart}_{f(i)-1}c_{i+1} \to \pi^{\heart}_{f(i)-1}c_{i}$ to be
	injective. This is exactly the condition that $\pi^{\heart}_{k}c_i = \pi^{\heart}_{k,i-k}c$ vanish when $f(i-1)\leq k<f(i)$ and $\pi^{\heart}_{f(i)-1}c_{i+1} = \pi^{\heart}_{f(i)-1,i-f(i)+2}c$ has no simple $\tau$-torsion.
\end{proof}

\begin{exm}\label{exm:slopeatstr}
	Let $f(i) = \ceil{a i}$ where $a \geq 0$. This gives rise to the \textit{slope $\frac {1-a}a$ $t$-structure}, whose truncation functors we denote $\tau^{/a}_{\geq0}, \tau^{/a}_{\leq 0}$.
\end{exm}

\begin{exm}\label{exm:slicentstr}
	Let $f(i) = 0$ for $i\leq 0$ and $f(i) = \ceil{\frac {i} 2}$ for $i >0$. This gives rise to the \textit{v $t$-structure}, whose truncation functors we denote $\tau^{v}_{\geq0},\tau^{v}_{\leq0}$.
\end{exm}

The slope $\frac{1-a}{a}$ and v $t$-structures satisfy the conditions of \Cref{lem:tsymmon} and \Cref{lem:ttruncate}, so are compatible with the symmetric monoidal structure, and can be computed by truncating level-wise. The reason for the name slope is that in the Adams grading, the homotopy groups of objects in the heart of this $t$-structure lie along a line of slope $\frac {1-a}{a}$. 
The v $t$-structure is named so because the curve it describes is the vanishing curve on the homotopy groups of the $\BP$-synthetic sphere at the prime $2$.

\begin{exm}\label{exm:filtstrs}
	We now specialize \Cref{exm:slopeatstr} to obtain two $t$-structures we use here.
	
	Taking $a = 0$, we get the \textit{constant} $t$-structure, whose connective cover functor $\tau^{\const}_{\geq0}$ just takes connective cover on each filtered piece.
	
	Taking $a = 1$, we get the \textit{diagonal} $t$-structure, whose connective cover functor $\tau^{d}_{\geq 0}$ is given by taking the $i^{th}$-connective cover on the $i^{\text{th}}$ filtered piece.
\end{exm}

The functor $(-)^{\const}:\mathcal C \to \Fil(\mathcal C)$ is the symmetric monoidal functor given by the constant filtered object.

\subsection{Filtrations on rings of interest}
We now specialize to the case $\mathcal C = \Sp$ with its standard symmetric monoidal structure. We begin by constructing $j_{\zeta}$ as a filtered ring. We use $\tau_{\geq*}(-)$ to denote the composite functor $\tau_{\geq 0}^{d}((-)^{\const})$. Indeed, $\tau_{\geq i}(-)$ is the $i^{\text{th}}$ filtered piece of this functor.

We now use $\tau_{\geq*}(-)$ to obtain a filtration on $\ell_p,j_{\zeta},\ko_2$, and $j$ for $p>2$. We use $R^{\fil}$ to denote these rings equipped with these filtrations, and $R^{\gr}$ to denote the associated graded algebras.

\begin{dfn}\label{defn:Zpfil}
	Let $\mathbb Z_p^{\fil}$ be the ring of $p$-adic integers with the $p$-adic filtration. It is a filtered $\mathbb E_\infty$-ring since it is a commutative ring in the heart of the constant $t$-structure. Its associated graded ring is $\mathbb F_p[v_0]$, where $v_0\in\pi_{0,1}\mathbb Z_p^{\gr}$. We write $\widetilde{v_0}\in\pi_{0,1}\mathbb Z_p^{\fil}$ for the class of filtration $1$ detecting $p\in\mathbb Z_p$, which projects to $v_0$ in the associated graded.
\end{dfn}

\begin{dfn}\label{defn:lkofilt}
	For $p>2$, consider $\ell_p$, viewed as an $\EE_{\infty}$-ring equipped with the $\ZZ$-action given by the Adams operation\footnote{The $\EE_{\infty}$-structure on Adams operations can be constructed from those on $\KU_p = \ZZ_p^{\times}$ by taking $\FF_p^{\times}$ fixed points and connective cover} $\Psi^{1+p}$, and for $p=2$, consider it with the $\ZZ\times C_2$-action given by the Adams operations $\Psi^{3}$, $\Psi^{-1}$.
	
	We now define most of our filtered $\EE_{\infty}$-rings of interest:
	\begin{itemize}
		\item $\ell_p^{\fil}:= \tau_{\geq*}\ell_p$
		\item $\ko_2^{\fil}:= \tau_{\geq0}^{v}((\ell_2^{\fil})^{hC_2})$ 
		\item $j_{\zeta,k}^{\fil}:= (\ell_p^{\fil})^{hp^k\ZZ}$ for $p>2$ and $(\ko_2^{\fil})^{h2^k\ZZ}$ for $p=2$
		\item $ju_{\zeta,k}^{\fil}:= (\ell_2^{\fil})^{h2^k\ZZ}$
		\item $j_{k}^{\fil}:=\tau_{\geq0}^{\const}(j_{\zeta,k}^{\fil})$ for $p>2$.
	\end{itemize}
	In the case $k=0$, we just write $j_{\zeta}^{\fil}, ju^{\fil}, j^{\fil}$, and we remove $\fil$ from the notation if we want to denote the underlying $\EE_{\infty}$-ring. For example, we write $j_{\zeta,k} = \ell_p^{hp^k\ZZ}$.
\end{dfn}

\begin{rmk}
	The filtrations of \Cref{defn:lkofilt} aren't as `fast' as they can possibly be. Namely, the spectra in the filtrations only change every multiple of $2p-2$ filtrations. Speeding up the filtration doesn't affect very much related to the filtration in any case.
\end{rmk}

\begin{rmk}
	For $p>2$, it is also possible to use variants of the \textit{Adams} filtration on the various rings of study, as in \cite[Section 4.3]{hahn2022redshift}, which would avoid the use of two filtrations. However this doesn't work as well at the prime $2$, since the Adams filtration is poorly suited to studying $\ko_2$'s $\THH$.
\end{rmk}
The key properties of these filtrations that we use is that the associated graded algebras mod $p$ are easy to describe.
\begin{lem}\label{lem:Zalg}
	The associated graded algebras of filtered rings defined in \Cref{defn:lkofilt} are $\EE_{\infty}$-$\ZZ$-algebras.
\end{lem}
\begin{proof}
	The $0$'th piece of every associated graded algebra is coconnective with $\pi_0 = \ZZ_p$, so the unit map from $\SP^{0,0}$ factors canonically through $\ZZ$, giving it a canonical $\EE_{\infty}$-$\ZZ$-algebra structure.
\end{proof}

\begin{lem}\label{lem:modp}
	For $p>2$, there are isomorphisms of graded $\EE_{\infty}$-$\FF_p$-algebras 
	\begin{align*}
	\ell_p^{\gr}/p &\cong \FF_p[v_1]\\
	j_{\zeta,k}^{\gr}/p&\cong \FF_p[v_1]\otimes_{\FF_p} \FF_p^{B\ZZ}
	\end{align*}
	and for $p=2$, there are isomorphisms of graded $\EE_{\infty}$-$\FF_2$-algebras
	\begin{align*}
	j_{\zeta,k}^{\gr}/2 &\cong (\ko_2^{\gr}/2)\otimes_{\FF_2} \FF_2^{B\ZZ}\\
	ju_{\zeta,k}^{\gr}/2 &\cong \FF_2[v_1]\otimes_{\FF_2} \FF_2^{B\ZZ}\\
	\ko_2^{\gr}/2 &\cong \tau_{\geq 0}^{v} (\FF_2^{BC_2}\otimes_{\FF_2}\FF_2[v_1]).
	\end{align*}
\end{lem}

\begin{figure}[hbt!]
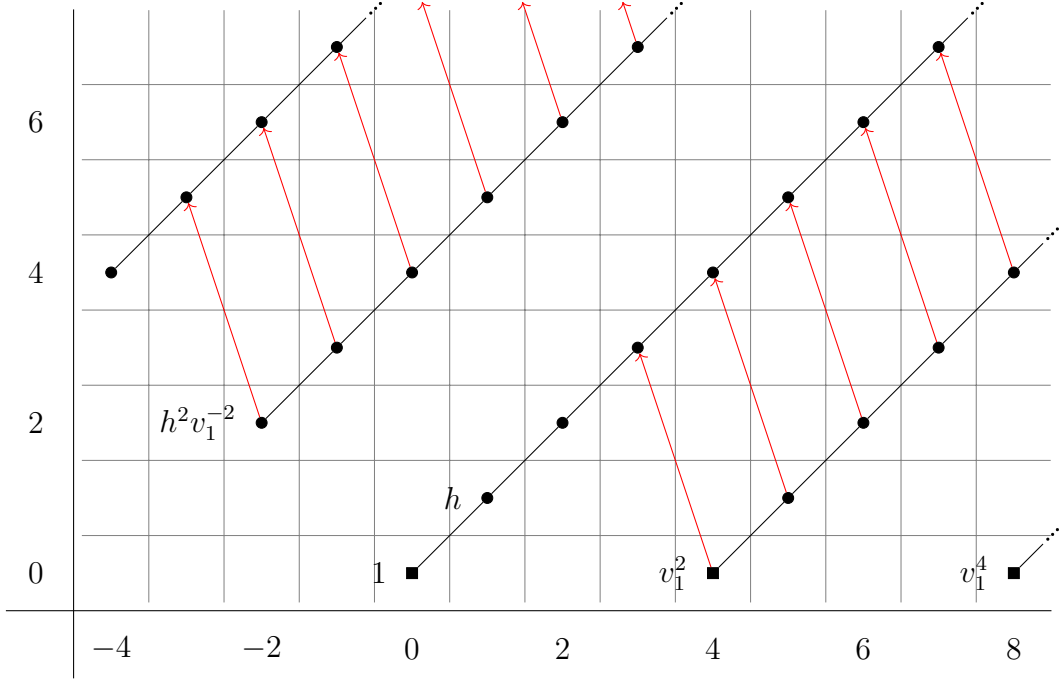

	\centering
	\scalebox{1.00}{
		\printpage[ name = kuhc2fil, page = 3 ]
	}
	\caption{
		Above is the $E_2$-page of the spectral sequence associated with the filtered ring $(\tau_{\geq*}\ell_2)^{hC_2}$, which embeds into the homotopy fixed point spectral sequence for $\KO_2$. The dots indicate a copy of $\FF_2$, the rectangles indicate a copy of $\ZZ_2$, and the diagonal lines indicate multiplication by $h$. The spectral sequence collapses at the $E_3$-page.
	}
	\label{fig:kuhc2fil}    
\end{figure}

\begin{proof}
	$\ell_{p}^{\gr}$ is the associated graded of the Postnikov filtration, which is $\ZZ_p[v_1]$, where the grading of $v_1$ is its topological degree, namely $2p-2$. Reducing mod $p$, we get the claim about $\ell^{\gr}/p$. The $\ZZ$-action on $\ell_{p}^{\gr}$ is the action of $\Psi^{1+p}$ on the homotopy of $\ell_p$. It is a ring automorphism sending $v_1$ to $(1+p)^{p-1}v_1$, which in particular is trivial modulo $p$. Since $\ell_p^{\gr}$ is a discrete object (it is in the heart of the diagonal $t$-structure), it follows that the action on $\ell_{p}^{\gr}/p$ is trivial, giving the claimed identification of $j_{\zeta,k}^{\gr}$ for $p>2$ and $ju_{\zeta,k}^{\gr}$ for $p=2$.
	
	For $p=2$, we first recall that in the homotopy fixed point spectral sequence for $\KO_2\cong \KU_2^{hC_2}$, all differentials are generated under the Leibniz rule by the differential $d_3v_1^2 = \eta^3$, where $\eta$ is represented by the class in $H^1(C_2;\pi_2\KU_2)$. The spectral sequence for $\ell_2^{hC_2} = \ku_2^{hC_2}$, displayed in \Cref{fig:kuhc2fil}, embeds into this, after a page shift. Thus, we see that everything in $\pi_{**}(\ku^{\fil})^{hC_2}$ above the line of slope $1$ intercept zero is either in negative underlying homotopy or doesn't have $\tau$-multiples on or below the line of slope $1$ intercept $2$. We learn that the bigraded homotopy ring of $(\ell^{\fil}_2)^{hC_2}$ is  
    \[
        \ZZ_2[x,\eta,\tau,b,v_1^4]/(b^2-4v_1^4,\eta^3\tau^2,2\eta,2x,x\eta\tau^2, v_1^4x-\eta^4,\eta b),
    \]
 where $x$ represents $v_1^{-4}\eta^4$, and $b$ represents $2v_1^2$.
	
	By applying \Cref{lem:grtruncate}, we learn that the connective cover $\tau_{\geq 0}^{v}$ can be computed the level of associated graded, and that this even holds after taking the cofiber by $2$. The $C_2$-action on $\ku_2^{\gr}/2$ is trivial, so indeed $\ko_2^{\gr}/2 \cong \gr(\tau_{\geq 0}^{v} (\FF_2^{BC_2}\otimes_{\FF_2}\FF_2[v_1]))$. For $j_{\zeta,k}^{\gr}/2$, we just observe that the residual $\ZZ$-action is also trivial.
\end{proof}
\begin{sseqdata}[ name = kuhc2fil, xscale=1.0, yscale=1.0, x range = {-4}{8}, y range = {0}{7}, x tick step = 2, y tick step = 2, class labels = {left}, classes = fill, grid = crossword, Adams grading, lax degree]
	\class["1",rectangle](0,0)
	\class["h"](1,1)\structline
	\class(2,2) \structline
	\class(3,3)\structline
	\class(4,4)\structline
	\class(5,5)\structline
	\class(6,6)\structline
	\class(7,7)\structline
	\class(8,8)\structline
	
	\class[rectangle,"v_1^2"](4,0)
	\class(5,1)\structline
	\class(6,2) \structline
	\class(7,3)\structline
	\class(8,4)\structline
	\class(9,5)\structline
	\class(10,6)\structline
	\class(11,7)\structline
	
	\class[rectangle,"v_1^4"](8,0)
	\class(9,1)\structline
	\class(10,2) \structline
	\class(11,3)\structline
	\class(12,4)\structline
	\class(13,5)\structline
	\class(14,6)\structline
	\class(15,7)\structline
	
	\class["h^2v_1^{-2}"](-2,2)
	\class(-1,3)\structline
	\class(0,4)\structline
	\class(1,5)\structline
	\class(2,6)\structline
	\class(3,7)\structline
	\class(4,8)\structline
	\class(5,9)\structline
	\class(6,10)\structline
	\class(7,11)\structline
	
	\class(-4,4)
	\class(-3,5)\structline
	\class(-2,6)\structline
	\class(-1,7)\structline
	\class(0,8)\structline
	\class(1,9)\structline
	\class(2,10)\structline
	\class(3,11)\structline
	
	\class(-6,6)
	\class(-5,7)\structline
	\class(-4,8)\structline
	\class(-3,9)\structline
	\class(-2,10)\structline
	\class(-1,11)\structline
	
    \d[red]3(4,0)(3,3)
    \d[red]3(5,1)(4,4)
    \d[red]3(6,2)(5,5)
    \d[red]3(7,3)(6,6)
    \d[red]3(8,4)(7,7)
    
    \d[red]3(-2,2)(-3,5)
    \d[red]3(-1,3)(-2,6)
    \d[red]3(0,4)(-1,7)
    \d[red]3(1,5)(0,8)
    \d[red]3(2,6)(1,9)   
    \d[red]3(3,7)(2,10)
	
\end{sseqdata}

\begin{sseqdata}[ name = juzeta2, xscale=1.0, yscale=1.0, x range = {-2}{11}, y range = {-8}{1}, x tick step = 2, y tick step = 2, class labels = {left}, grid = crossword, Adams grading, lax degree]

	\class["\lambda_1'",rectangle](3,-3)
	\class[rectangle](2,-2)
	\class[rectangle](7,-7)
	\class[rectangle](6,-6)
	\class[rectangle](11,-11)
	\class[rectangle](10,-10)
	
	\class[rectangle](6,-4)
	\class[rectangle](5,-3)
	\class[rectangle](10,-8)
	\class[rectangle](9,-7)
	
	\class["1",rectangle](0,0)
	\class["\zeta",rectangle](-1,1)
	\class["\sigma^2v_1",rectangle](4,-4)
	\class[name = nl,rectangle](3,-3)
	\class[rectangle](8,-8)
	\class[rectangle](7,-7)
	
	\class["dv_1",rectangle](3,-1)
	\class[rectangle](2,0)
	\class[rectangle](7,-5)
	\class[rectangle](6,-4)
	\class[rectangle](11,-9)
	\class[name = s,rectangle](10,-8)
	
	\d[red]2(11,-11)(s)

	\d[red]2(4,-4)(3,-1)
	\d[red]2(nl)(2,0)
	\d[red]2(6,-6)(5,-3)
	\d[red]2(7,-7)(6,-4)
	\d[red]2(7,-7)(6,-4)
		\d[red]2(7,-7)(6,-4)

	
\end{sseqdata}

\begin{sseqdata}[ name = thhko, xscale=1.0, yscale=1.0, x range = {-1}{11}, y range = {-8}{1}, x tick step = 2, y tick step = 2, class labels = {left}, fill, grid = crossword, Adams grading, lax degree]

	\class["1"](0,0)
	\class["\sigma^2v_1"](4,-4)
	\class(8,-8)
	
	\class["dv_1"](3,-1)
	\class(7,-5)
	\class(11,-9)
	
	\class["d\eta"](2,0)
	\class[name = b](6,-4)
	\class[name = d](10,-8)
	
	\class(5,-1)
	\class(9,-5)
	\class(13,-9)
	
	\class["\sigma^2\eta"](3,-3)
	\class(7,-7)
	\class(11,-11)
	
	\class[name = a](6,-4)
	\class[name = c](10,-8)
	\class(14,-12)
	
	\class(5,-3)
	\class(9,-7)
	\class(13,-11)
	
	\class(8,-4)
	\class(12,-8)
	\class(16,-12)
	\d[red]2(3,-3)(2,0)
	\d[red]2(4,-4)(3,-1)
	\d[red]2(9,-7)(8,-4)
	\d[red]2(7,-7)(b)
	\d[red]2(7,-7)(a)
	\d[red]2(a)(5,-1)
	\d[red]2(b)(5,-1)
	
	\d[red]2(c)(9,-5)
	\d[red]2(11,-11)(d)

	
\end{sseqdata}

\begin{sseqdata}[ name = thhj, xscale=.333, yscale=.333, x range = {-2}{40}, y range = {-24}{1}, x grid step = 4, y grid step = 3, x tick step = 8, y tick step = 6, class labels = {left}, fill, grid = crossword, Adams grading, lax degree]

	\class["1"](0,0)
	\class["\alpha_1"](3,1)
	\class["\gamma_p(d\alpha_1)"](12,0)
	\class(15,1)
	\class(24,0)
	\class(27,1)
	\class["\gamma_{p^2}(d\alpha_1)"](36,0)
	\class(39,1)
	
	\class["a"](13,-5)
	\class(16,-4)
	\class(25,-5)
	\class(28,-4)
	\class(37,-5)
	\class(40,-4)
	\class(49,-5)
	\class(52,-4)
	
	\class["\lambda_2"](17,-13)
	\class(20,-12)
	\class(29,-13)
	\class(32,-12)
	\class(41,-13)
	\class(44,-12)
	\class(53,-13)
	\class(56,-12)
	
	\class(30,-18)
	\class(33,-17)
	\class(42,-18)
	\class(45,-17)
	\class(54,-18)
	\class(57,-17)
	\class(66,-18)
	\class(69,-17)
	
	\class["(\sigma^2v_1)^p"](18,-18)
	\class(21,-17)
	\class(30,-18)
	\class(33,-17)
	\class(42,-18)
	\class(45,-17)
	\class(54,-18)
	\class(57,-17)
	
	\class(31,-23)
	\class(34,-22)
	\class(43,-23)
	\class(46,-22)
	\class(55,-23)
	\class(58,-22)
	\class(67,-23)
	\class(70,-22)
	
	\class(35,-31)
	\class(38,-30)
	\class(47,-31)
	\class(50,-30)
	\class(59,-31)
	\class(62,-30)
	\class(51,-31)
	\class(54,-30)
	
	\d[purple]2(17,-13)(16,-4)
	\d[purple]2(29,-13)(28,-4)
	\d[purple]2(41,-13)(40,-4)
	
	\d[purple]2(35,-31)(34,-22)
	\d[purple]2(29,-13)(28,-4)
	\d[purple]2(41,-13)(40,-4)
	
\end{sseqdata}

\begin{rmk}
	At the prime $2$, it is possible to define $j$ as a filtered $\EE_{\infty}$-ring, but we do not study this in this paper. One can define its underlying $\EE_{\infty}$-ring as the pullback
	
	\begin{center}
		\begin{tikzcd}
			j \pullback \ar[r]\ar[d] &\ko_2^{h\ZZ} \ar[d]\\
			\tau_{\leq2}\SP_2 \ar[r] & (\tau_{\leq2}\ko_2)^{h\ZZ}
		\end{tikzcd}
	\end{center}
	and then consider the underlying filtered $\EE_{\infty}$-ring of $\nu_{BP}(j)$ where $\nu_{BP}$ is the synthetic analogue functor of \cite{pstrkagowski2022synthetic}.
\end{rmk}

Finally, we show convergence properties of our $\THH$ applied to the filtrations we use.	Given a filtered spectrum $X\in\Fil(\Sp)$, the spectral sequence associated with $X$ converges conditionally if and only if $\lim_i X_i = 0$. This is equivalent to asking that $X$ is $\tau$-complete, where $\tau$ is in $\pi_{0,-1}\SP^{0,0}$.

The following lemma shows completeness for $\THH$ with respect to all of the filtrations constructed in this section.
\begin{lem}\label{lem:postnikovcomplete}
	Suppose that $R$ is a filtered ring such that the $i$-th filtered piece $R_i$ is $(-1+ci)$-connective for every $i$ and some fixed $c>0$. Then, the $i$-th filtered piece of $\THH(R)$ is also $(-1+ci)$-connective, so in particular the filtration on $\THH(R)$ is complete.
\end{lem}

\begin{proof}
	Note that $\overline R=\cof(\mathbb S^{0,0}\to R)$ satisfies the same conditions of the statement. The filtration from the cyclic bar construction (see \cite[III.2]{nikolaus2018topological}) gives us an increasing filtration on $\THH(R)$ with $k$-th associated graded piece $\Sigma^k R\otimes \overline{R}^{\otimes k}$. The $i$-th filtered piece of $\Sigma^k R\otimes \overline{R}^{\otimes k}$ is $(-1+ci)$-connective since it is a colimit of spectra of the form $\Sigma^k R_{j_0}\otimes \overline{R}_{j_1}\otimes\cdots\otimes \overline{R}_{j_k}$ with $j_0+\cdots+j_k\geq i$, which has connectivity of at least
	\[
	k + \sum_{s=0}^k (-1 + cj_s) \geq -1 + ci. \qedhere
	\]
\end{proof}

\begin{lem}\label{lem:padiccomplete}
	Let $R$ be a (possibly graded) $\mathbb E_1$-$\mathbb Z_p$-algebra. Then, the filtration on the filtered ring
	\[
		\THH(R\otimes_{\mathbb Z_p}\mathbb Z_p^{\fil})/\widetilde{v_0}.
	\]
	is complete and its associated graded ring is concentrated in two filtration degrees $t=0,1$. Informally, the filtration is of the form
	\[
		\cdots\to0\to0\to I\to \THH(R)/p
	\]
	for some (possibly graded) spectrum $I$. In particular, the associated spectral sequence collapses at the $E_2$-page.
\end{lem}
\begin{proof}
	By using the symmetric monoidality of $\THH$ and the fact that $p=0$ in $\ZZ_p^{\fil}/\tilde{v}_0$, we obtain an equivalence
	\[
	\THH(R\otimes_{\mathbb Z_p}\mathbb Z_p^{\fil})/\tilde{v}_0 \cong (\THH(R)/p) \otimes_{(\THH(\mathbb Z_p)/p)}\THH(\mathbb Z_p^{\fil})/\tilde{v}_0.
	\]
	Since the conclusion of the statement is stable under base-change along trivially filtered rings, the statement reduces to the case $R=\mathbb Z_p$.
	
	For $R=\mathbb Z_p$, the associated graded is $\THH(\FF_p[v_0])/v_0$, which has homotopy ring $\FF_p[\sigma^2p]\otimes \Lambda[dv_0]$ (see \Cref{exm:THHZ}), which is indeed in filtrations $\leq1$. It remains to see that $\THH(\ZZ_p^{\fil})/\tilde{v}_0 = \THH(\ZZ_p^{\fil};\FF_p)$ has a complete filtration. It suffices to show that $\THH(\ZZ_p^{\fil};\FF_p)\otimes_{\THH(\ZZ_p)}\ZZ_p \cong \THH(\ZZ_p^{\fil}/\ZZ_p;\FF_p)$ has a complete filtration, since $\THH(\ZZ_p)$ is built from $\ZZ_p$ via extensions and limits that are finite in each degree, and completeness of the filtration can be checked degreewise. The $n$th associated graded term of the cyclic bar construction computing this is $$\Sigma^n (\ZZ_p^{\fil})^{\otimes_{\ZZ_p} n}\otimes_{\ZZ_p}\FF_p \cong  \Sigma^n(\ZZ_p^{\fil}\otimes_{\ZZ_p}\FF_p)^{\otimes_{\FF_p} n}$$
	
	$\ZZ_p^{\fil}\otimes_{\ZZ_p}\FF_p$ is complete since it is $\FF_p$ in each nonnegative degree, with transition maps $0$, or in other words, it is a direct sum $\FF_p\oplus \bigoplus_{1}^{\infty}\Sigma^{0,i}\FF_p/\tau$. It follows that its tensor powers over $\FF_p$ are also sums of $\FF_p$ in each degree with transition maps $0$ in positive filtration, so are complete. Since only finitely many terms in the cyclic bar complex contribute to each degree of $\THH$, we learn that the $\THH$ is complete.
\end{proof}

\section{Tools for understanding $\THH$}\label{sec:tools}

In this section, we explain some general tools which we use in understanding $\THH$.

\subsection{Suspension operation in THH}
We begin by reviewing and proving some basic facts about the suspension maps, which are studied in \cite[Section A]{hahn2022redshift}. Let $R$ be an $\mathbb E_1$-algebra in a presentably symmetric monoidal stable category $\mathcal C$. By \cite[Section A]{hahn2022redshift}, there are natural maps
\begin{align}
	\sigma:\Sigma\fib(1_R)&\to R\otimes R\label{eq:sigmas0}\\
	\sigma^2:\Sigma^2\fib(1_R)&\to \THH(R)\label{eq:sigmas1}
\end{align}
where $1_R$ is the unit map of $R$. Note that the first map is defined by the diagram
\[
\begin{tikzcd}[column sep = huge]
	\unit\ar[r]\ar[d,"1_R"]&0\ar[d]\\
	R\ar[r,"\id\otimes1_R-1_R\otimes\id"]&R\otimes R
\end{tikzcd}
\]
and that it factors through $\fib(\mu)\to R\otimes R$ where $\mu:R\otimes R\to R$ is the multiplication map.

Let $I$ be an object of $\mathcal C$ with a map $I\to R\otimes R$ and nullhomotopies of the composites
\begin{align*}
	I&\to R\otimes R\xrightarrow\mu R\\
	I&\to R\otimes R\xrightarrow{\mu\circ T}R,
\end{align*}
where $T:R\otimes R\to R\otimes R$ is the exchange map. Then, we obtain a map
\[
	\Sigma I\to \THH(R)
\]
by the commutative diagram
\begin{equation}
\begin{tikzcd}
	I\ar[rr]\ar[rd]\ar[dd]&&0\ar[d]\\
	&R\otimes R\ar[r,"\mu"]\ar[d,"\mu\circ T"]&R\ar[d,"1\otimes\id"]\\
	0\ar[r]&R\ar[r,"\id\otimes1"]&R\otimes_{R\otimes R^{op}}R.
\end{tikzcd}\label{eq:cdsusp}
\end{equation}
By the proof of \cite[Lemma A.3.2]{hahn2022redshift}, if $I=\Sigma\fib(1_R)$ and the map $I\to R\otimes R$ is given by \eqref{eq:sigmas0}, then the induced map $\Sigma I\to\THH(R)$ is the map \eqref{eq:sigmas1}.

\begin{dfn}
Let $X$ be a spectrum. Given a class $x\in\pi_\ast(X\otimes1)$ and a lift $\widetilde{x}\in\pi_\ast(X\otimes \fib(1_R))$ we shall write $\sigma x\in\pi_{\ast+1}(X\otimes R\otimes R)$ and $\sigma^2 x\in\pi_{\ast+2}(X\otimes\THH(R))$ for the image of $\widetilde{x}$ under the maps \eqref{eq:sigmas0} and \eqref{eq:sigmas1}. The notation is ambiguous since we need to choose a lift $\widetilde x$, but these lifts will often be well-defined.

We shall write $d$ for
\[
\pi_\ast(X\otimes R)\to\pi_{\ast+1}(X\otimes \THH(R))
\]
induced by the map of spectra $\Sigma R\to \Sigma^2\fib(1_R)\to \THH(R)$.
\end{dfn}

\begin{rmk}\label{rmk:fibmul}
If $R$ is homotopy commutative in addition to being an $\mathbb E_1$-algebra, then we can set $I=\fib(\mu)$ in \eqref{eq:cdsusp} and obtain a map
\begin{equation}
	\sigma:\Sigma\fib(\mu)\to \THH(R),\label{eq:sigma}
\end{equation}
which is functorial on $R$ and the homotopy\footnote{The same construction is studied in \cite[Variant A.2.2]{hahn2022redshift}, but we believe that additional hypotheses are required to make sense of their argument. For example, $R$ is only assumed to be an $\mathbb E_1$-ring in their generality, but an assumption such as homotopy commutativity of $R$ is needed to ensure that the composite
	\[
	\fib(\mu)\to R\otimes R\xrightarrow{\mu\circ T}R
	\]
	is nullhomotopic. In their notation, we would need to assume, for example, that there is a homotopy $1_k\simeq 1_k^\tau$. This does not affect any other part of their work since they only use rings that have enough structure.} $\mu\simeq \mu\circ T$. Then, the map \eqref{eq:sigmas1} is the composite
\[
\Sigma^2\fib(1_R)\xrightarrow{} \Sigma\fib(\mu)\xrightarrow{\sigma}\THH(R)
\]
Here the first map is obtained from using the null homotopy of the map (1) composed with the multiplication.

If $X$ is a spectrum, given a class $y\in\pi_\ast(X\otimes R\otimes R)$ and a lift $\widetilde{y}\in\pi_\ast(X\otimes \fib(\mu))$, we shall write $\sigma y\in \pi_{\ast+1}(X\otimes \THH(R))$ for the image of $\widetilde y$ under the map \eqref{eq:sigma}. Then, we have $dx= \sigma((\eta_L-\eta_R)x)$ for $x\in\pi_\ast(X\otimes R)$, where $\eta_L$ and $\eta_R$ are the left and right units of $R\otimes R$, respectively.
\end{rmk}

\begin{lem}[{\cite[Prop. 5.10]{angeltveit2005hopf}}]\label{lem:dLeibniz}
	Let $X$ be a homotopy unital ring spectrum and $R$ be an $\mathbb E_2$-algebra in $\mathcal C$. Then, $d$ satisfies the Leibniz rule
	\[
		d(xy) = d(x)y + (-1)^{|x|}xd(y)
	\]
	for any $x,y\in \pi_\ast(X\otimes R)$.
\end{lem}
\begin{proof}
	By \cite[Example A.2.4]{hahn2022redshift}, the map $d$ can be identified with the map
	\[
		S^1_+\otimes R \to \THH(R)
	\]
	induced by the unit map $R\to\THH(R)$ and the $S^1$-action on $\THH(R)$. Since the map $R\to\THH(R)$ is a map of $\mathbb E_1$-rings, the $S^1$-action on the target gives an $S^1$-family of ring maps, and so we obtain a map of $\mathbb E_1$-rings
	\begin{equation}\label{eq:dcircle}
		R\to \lim_{S^1}\THH(R) = DS_+^1\otimes\THH(R) = \THH(R) \oplus \Sigma^{-1}\THH(R)
	\end{equation}
	given by the sum of the identity map and $d$. Here, $DS_+^1$ is the Spanier-Whitehead dual of $S^1$ with the algebra structure given by the diagonal map of $S^1$.
	
	The homotopy ring of $DS_+^1$ is given by
	\[
		\pi_\ast(DS_+^1) = (\pi_\ast S^0)[t]/(t^2)
	\]
	with $|t|=-1$. Since \eqref{eq:dcircle} is a ring map, taking the $X$-homology, we have
	\[
		1\otimes xy + t\otimes d(xy)=(1\otimes x + t\otimes dx )(1\otimes y+t\otimes dy)
	\]
	for $x,y\in \pi_\ast (X\otimes R)$. Expanding it using $t^2=0$ gives us the desired Leibniz rule.
\end{proof}

Our use of the symbol $d$ recovers the use in the HKR theorem. Recall that a strict Picard element $\sL$ of a symmetric monoidal category $\mathcal C$ is a map of spectra $\ZZ \to \pic(\mathcal C)$. Given such a strict Picard element, viewing it as a symmetric monoidal functor $\ZZ \to \mathcal C$, the colimit of the composite
\[
\mathbb N\to\mathbb Z\to\mathcal C
\]
is an $\EE_{\infty}$-algebra in $\mathcal C$ which we denote $\unit[x]$, where $x$ is a class in the Picard graded homotopy in the degree of $\sL$.

\begin{lem}[HKR isomorphism]\label{lem:hkrdx}
	Let $C$ be a presentably symmetric monoidal stable category with a strict Picard element $\sL$.
	Let $\unit[x]$ denote the polynomial algebra on a class $x$ in degree $\sL$. Then $\HoH(\unit[x])$ is a free $\unit[x]$-module on $1$ and $dx$.
\end{lem}

\begin{proof}
	The universal example of such a $C$ is graded spectra, where $\unit[x]$ is the graded polynomial algebra $\Sigma^{\infty}_+\NN$, so it suffices to prove it there. But now this follows from from the Kunneth spectral sequence computing $\pi_*\THH(\SP[x]) = \pi_*\SP[x]\otimes_{\SP[x_1,x_2]}\SP[x]$, since $dx$ is $\sigma((\eta_L-\eta_R)(x))$.
\end{proof}

We now explain some basic $\THH$ computations involving the suspension map.

\begin{exm}[B\"okstedt periodicity]
	The fundamental computation of B\"okstedt states that the ring $\pi_\ast\THH(\mathbb F_p)$ is isomorphic to $\mathbb F_p[\sigma^2p]$.
\end{exm}

\begin{lem}\label{lem:diffonsigma}
	Let $R\in\Fil(\Sp)$ be a filtered $\EE_1$-ring and $X\in\Sp$ a spectrum. Let $y\in \pi_{k,r-k}(R\otimes X), \; x\in \pi_k X$ be classes such that $\tau^r y =x\in\pi_{k,-k}(R\otimes X)$. 

	Then there is a choice of nullhomotopy of $x$ in $\THH(\gr R)\otimes X$ such that in the spectral sequence for $\THH(R)\otimes X$, the corresponding element $\sigma^2x$ on the $E_1$-page survives to the $E_r$-page and has $d_r$-differential $d_r(\sigma^2 x)=\pm d y$.
\end{lem}
\begin{proof}
	A choice of homotopy $\tau^r y \sim x$ in $R\otimes X$ becomes in  $\cof(\SP^{0,0}\to R)\otimes X$ a choice of nullhomotopy of the image of $\tau^r y$, which corresponds to a map $\Sigma^{|y|}\cof(\tau^r) \to \cof(\SP^{0,0}\to R)\otimes X$. This map of filtered spectra gives a map of the associated spectral sequences, and in the spectral sequence for $\cof(\tau^r)$, there is a $d_r$-differential between the two spheres on the associated graded.
	
	We claim the image of the two shifts of $\cof\tau$ in the map $$\Sigma^{|y|}(\cof(\tau) \oplus \Sigma^{1,-({r+1})}\cof (\tau)) \cong \Sigma^{|y|}\cof(\tau^r)\otimes \cof(\tau) \to \cof(\SP^{0,0}\to R)\otimes X\otimes \cof(\tau)$$ correspond to the image of $y$ and the suspension of a nullhomotopy of $x$ under the map $\SP^{0,0} \to \gr R$.
	
	The claim that the first $\cof\tau$ is sent to $y$ is clear by construction, and the claim that the second $\cof\tau$ is sent to the suspension of a nullhomotopy of $x$ follows since on associated graded our original homotopy $\tau^ry\sim x$ becomes a nullhomotopy of $x$.
	
	It then follows that there is a $d_r$ differential between these two classes.
	
	Composing with the filtered map $$\Sigma\cof(\SP^{0,0} \to R)\otimes X\cong \Sigma^2 \fib(\SP^{0,0} \to R)\otimes X \xrightarrow{\sigma^2} \THH(R)\otimes X$$ of \Cref{eq:sigmas1}, $y$ gets sent to $dy$ and the nullhomotopy of $x$ gets sent to $\sigma^2x$ (up to a possible sign), giving the desired differential in the spectral sequence for $\THH(R)\otimes X$.
    
	Therefore, it is enough to prove that the connecting map sends $\widetilde x$ to $y$, and since the map $\pi_\ast(Z\otimes X_1)\to\pi_\ast(Z\otimes X_0)$ is injective, it is enough to prove that $\widetilde x$ is sent to $\eta_\ast(x)$ by the composite $F\to X_1\to X_0$. This composite is homotopic to $F\to \mathbb S\to X_0$ since the connecting map $F\to X_1$ is given by the nullhomotopy.
\end{proof}

\subsection{THH in the stable range}

Throughout this subsection, let $S$ be a connective $\EE_{\infty}$-algebra and $R$ be a connective $\mathbb E_1$-$S$-algebra.

In this section, we show that in the situation that the unit map $S \to R$ is highly connective, $\THH(R/S)$ in low degrees becomes relatively straightforward to understand. This is used later in \Cref{sec:j} to understand $\THH(j)$. Let $\Delta_n$ denote the subcategory of $\Delta$ consisting of ordinals of size $\leq n$.

\begin{lem}\label{prop:cycbarvanline}
	If the unit map $S\to R$ is $i$-connective, then the natural map 
	\[
	\colim_{\Delta^{op}_n}R^{\otimes_S{*+1}} \to \colim_{\Delta^{op}}R^{\otimes_S{*+1}} \cong \THH(R/S)
	\]
	is $(n+1)(i+2)-1$-connective.
\end{lem}
\begin{proof}
	Let $\overline{R}=\cof(S\to R)$ be the cofiber of the unit map. The $m^{th}$ term of the associated graded of the filtration coming from the cyclic bar construction is $\Sigma^m R\otimes_S \overline{R}^{\otimes_S{m}}$, which is $m(i+2)$-connective because $R$ is connective and $\overline{R}$ is $(i+1)$-connective. It follows that the cofiber of the map in question has an increasing filtration whose associated graded pieces are $m(i+2)$-connective for $m >n$. This implies the result.
\end{proof}

The above lemma gives a simple description of $\THH$ in low degrees.

\begin{prop}\label{prop:stableTHH}
	If the unit map $S \to R$ is $i$-connective, then the map
	\[
	\Sigma^2\fib(1_R) \oplus R \xrightarrow{\sigma^2\oplus 1} \THH(R/S)
	\]
	is $(2i+2)$-connective, where $\sigma^2$ is defined as in \eqref{eq:sigmas1}.
\end{prop}
\begin{proof}
	Consider the case $n=1$ in Lemma \ref{prop:cycbarvanline}. Then, we have an equivalence
	\[
	\colim_{\Delta^{op}_1} R^{\otimes_S{\ast+1}} \simeq \colim\left(
	\begin{tikzcd}
		R\otimes_S R\ar[r,"\mu"]\ar[d,"\mu\circ T"]&R\\
		R&
	\end{tikzcd}
	\right)
	\]
	(see \cite[Theorem 9.4.4]{munson2015cubical}), where $T$ is the exchange map, and this colimit maps into $\THH(R/S)$ by a $(2i+3)$-connective map.
	
	Therefore, it is enough to prove that the map
	\[
		\colim\left(
			\begin{tikzcd}
				\Sigma\fib(1_R)\oplus R\ar[r,"\mathrm{proj}_2"]\ar[d,"\mathrm{proj}_2"]&R\\
				R&
			\end{tikzcd}
			\right)
		\to
		\colim\left(
			\begin{tikzcd}
				R\otimes_S R\ar[r,"\mu"]\ar[d,"\mu\circ T"]&R\\
				R&
			\end{tikzcd}
			\right)
	\]
	is $(2i+2)$-connective, where the map $\Sigma\fib(1_R)\oplus R \to R\otimes_SR$ is $\sigma\oplus(1_R\otimes\mathrm{id})$ and the two maps $R\to R$ are the identities. The fiber of this map is
	\[
	\Sigma\fib(\Sigma\fib(1_R)\oplus R\xrightarrow{\sigma\oplus1} R\otimes_SR)
	\]
	which is $(2i+2)$-connective by the next lemma.
\end{proof}

\begin{lem}
	If the unit map $1_R:S\to R$ is $i$-connective, then the map
	\[
	\Sigma\fib(1_R)\oplus R\xrightarrow{\sigma\oplus1} R\otimes_SR
	\]
	is $(2i+1)$-connective.
\end{lem}
\begin{proof}
	This is equivalent to asking that the total cofiber of the following diagram
	\[
	\begin{tikzcd}
		S\otimes_SS\ar[r]\ar[d]&S\otimes_SR\ar[d]\\
		R\otimes_SS\ar[r]&R\otimes_SR
	\end{tikzcd}
	\]
	is $(2i+2)$-connective. This follows from the assumption since the total cofiber is $\Sigma^2 \fib(1_R)\otimes_{S}\fib(1_R)$, which is $(2i+2)$-connective since $\fib(1_R)$ is $i$-connective.
\end{proof}

\begin{cor}\label{cor:lambda1}
	The group $\pi_{2p-1}\THH(\mathbb Z_p)$ is isomorphic to $\mathbb Z/p$ and is generated by $\sigma^2\alpha_1$.
\end{cor}
\begin{proof}
	Since $\mathbb S_p\to\mathbb Z_p$ is $(2p-3)$-connective, the result follows from \Cref{prop:stableTHH}, which implies that $\sigma^2$ induces an isomorphism
	\[
		\mathbb Z/p=\pi_{2p-3}\fib(\mathbb S_p\to\mathbb Z_p)\simeq\pi_{2p-1}\THH(\mathbb Z_p).\qedhere
	\]
\end{proof}

\begin{cor}\label{cor:jconn}
	For $p>2$, the map
	\[
		j\oplus \Sigma^2\fib(\mathbb{S}_p\to j) \xrightarrow{1\oplus \sigma^2}\THH(j)
	\]
	is $(4p^2-4p-2)$-connective.
\end{cor}
\begin{proof}
For $p>2$, $\SP_p \to j$ is $2p^2-2p-2$-connective. This is because the first element of the fiber is $\beta_1$ (see for example \cite[Theorem 4.4.20]{ravenelgreen}) which is in that degree.
\end{proof}

\section{The $\THH$ of $j_{\zeta}$}

In this section, we compute $\THH(j_{\zeta})/(p,v_1)$ using the filtration constructed in \Cref{sec:filtration}. Let us first assume that $p$ is an odd prime. We shall discuss the case $p=2$ later in the section.

\subsection{THH of $\mathbb Z_p$ and $\ell_p$}
Before computing the $\THH$ of $j_\zeta$, we shall compute the $\THH$ of $\mathbb  Z_p$ modulo $p$ and the $\THH$ of $\ell_p$ modulo $(p,v_1)$ in this section, as a warm-up. They will be computed using the spectral sequences associated with $\THH(\mathbb Z_p^{\fil})$ and $\THH(\ell_p^{\fil})$. Later, we show that the computation of the spectral sequence for $\THH(j_{\zeta}^{\fil})$ looks the same. We note that the computations for $\ZZ_p$ and $\ell_p$ are well-known (see for example \cite[Theorem 5.12]{angeltveit2005hopf}).

\begin{lem}\label{lem:THHpoly}
	Let $k$ be a discrete commutative ring and let $R$ be a $\mathbb Z^m$-graded $\mathbb E_2$-$k$-algebra such that the homotopy groups of $R$ form a polynomial $k$-algebra
	\[
		\pi_\ast R = k[x_1,\dots,x_n]
	\]
	on even degree generators $x_1,\dots,x_n$. Then, there is an equivalence of $\mathbb Z^m$-graded $\mathbb E_1$-$\THH(k)$-algebras
	\[
		\THH(R)\simeq \THH(k)\otimes_k \HoH(k[x_1,\dots,x_n]/k).
	\]
\end{lem}
\begin{proof}
	Let $\mathbb S[x_1,\dots,x_n]$ be the $\mathbb Z^m$-graded $\mathbb E_2$-ring spectrum of \cite{rotation}. Then, by \cite[Prop. 4.2.1]{hahn2022redshift}, there is an equivalence of $\mathbb Z^m$-graded $\mathbb E_2$-$k$-algebras
	\[
		R\simeq k\otimes\mathbb S[x_1,\dots,x_n].
	\]
	Therefore, since $\THH$ is a symmetric monoidal functor $\Alg(\Sp)\to\Sp$, there is an equivalence of $\mathbb Z^m$-graded $\mathbb E_1$-$k$-algebras
	\[
		\THH(R) \simeq \THH(k) \otimes\THH(\mathbb S[x_1,\dots,x_n]),
	\]
	and the statement follows by base changing the second tensor factor on the right hand side along $\mathbb S\to k$.
\end{proof}

\begin{exm}\label{exm:THHZ}
	Consider the filtered spectrum $\THH(\ZZ_p^{\fil})/\tilde{v_0}$. Its associated graded spectrum is $\THH(\FF_p[v_0])/v_0$ and its underlying spectrum is $\THH(\ZZ_p)/p$. The $E_1$-page of the associated spectral sequence is $\mathbb F_p[\sigma^2p]\otimes\Lambda[dv_0]$ by \Cref{lem:THHpoly}. Note that $\sigma^2 p$ and $dv_0$ are in filtrations $0$ and $1$, respectively.

	By \Cref{lem:diffonsigma}, we have a differential $d_1(\sigma^2p)\doteq dv_0$ in the spectral sequence associated with the filtered ring $\THH(\mathbb Z_p^{\fil})$. Then, mapping to $\THH(\mathbb Z_p^{\fil})/\widetilde{v_0}$ and using the Leibniz rule, we can determine all differentials, and the $E_2$-page is isomorphic to $\mathbb F_p[(\sigma^2p)^p]\otimes\Lambda[(\sigma^2p)^{p-1}dv_0]$. There are no differentials in later pages by \Cref{lem:padiccomplete}.

	Therefore, the homotopy ring $\pi_\ast\THH(\mathbb Z_p)/p$ is isomorphic to $\mathbb F_p[\mu]\otimes\Lambda[\lambda_1]$ with $|\mu|=2p$ and $|\lambda_1|=2p-1$. By \cite[Propsition 6.1.6]{hahn2022redshift}, $\mu$ can be identified with $\sigma^2 v_1$\footnote{$v_1$ is not well defined at the prime $2$, but still exists: it is just not a self map of $\cof(2)$. It is generally defined as any element of $\pi_{2p-2}\SP/p$ whose $\BP$-Hurewicz image is $v_1$.}, where $v_1 \in \pi_{2p-2}\SP_p$ and $\lambda_1$ can be identified with $\sigma t_1$, in the sense of \Cref{rmk:fibmul}, where $t_1\in\pi_\ast(\mathbb Z\otimes\mathbb Z)$ is the image of $t_1\in\pi_\ast(\BP\otimes\BP)$ under the map $\BP\to\mathbb Z$. By \Cref{cor:lambda1}, we have $\lambda_1 \doteq \sigma^2\alpha_1$\footnote{Alternatively, if one knows that the $p$-Bockstein on $\mu$ is $\doteq \lambda_1$, one learns that $\sigma^2\alpha \doteq \lambda_1$ from the fact that the $p$-Bockstein on $v_1$ is $\alpha_1$ and the fact that $\sigma^2$ is compatible with the $p$-Bockstein (since it comes from a map of spectra).}.
\end{exm}

\begin{exm}\label{exm:THHell}
	Consider the filtered spectrum $\THH(\ell_p^{\fil})/(p,\widetilde{v_1})$, where $\widetilde{v_1}\in\pi_\ast\ell_p$ is the class of filtration $(2p-2)$. Its associated graded spectrum is $\THH(\mathbb Z[v_1])/(p,v_1)$ and its underlying spectrum is $\THH(\ell_p)/(p,v_1)$. By \Cref{lem:THHpoly}, the $E_1$-page of the associated spectral sequence is $\mathbb F_p[\sigma^2v_1]\otimes\Lambda[\lambda_1,dv_1]$. Note that the for degree reasons, the first and last page a differential can happen is the $E_{2p-2}$-page.
	
	Applying Lemma \ref{lem:diffonsigma}, there is a differential $d_{2p-2}\sigma^2 v_1 \doteq dv_1$ in the spectral sequence associated with the filtered spectrum $\THH(\ell_p^{\fil})/p$. Mapping to $\THH(\ell_p^{\fil})/(p,\widetilde{v_1})$ and using the Leibniz rule, we can determine the $d_{2p-2}$-differentials on powers of $\sigma^2 v_1$. The class $\lambda_1$ is a permanent cycle for degree reasons. Therefore, the $E_{2p-1}$-page is isomorphic to $\mathbb F_p[(\sigma^2 v_1)^p]\otimes\Lambda[\lambda_1, (\sigma^2v_1)^{p-1}dv_1]$. The classes $(\sigma^2v_1)^p,(\sigma^2v_1)^{p-1}dv_1$ are permanent cycles for degree reasons, so the spectral sequence degenerates at the $E_{2p-1}$-page.

	We let $\lambda_2$ denote a class detecting $(\sigma^2v_1)^{p-1}dv_1$, and $\mu$ denote a class detecting $(\sigma^2v_1)^p$.
	To check that there are no multiplicative extensions, we need to check $\lambda_1^2=\lambda_2^2=0$, which follows for degree reasons. The homotopy ring $\pi_\ast\THH(\ell_p)/(p,v_1)$ is thus isomorphic to $\mathbb F_p[\mu_1]\otimes\Lambda[\lambda_1,\lambda_2]$ where $\lambda_1$ and $\lambda_2$ can be identified with $\sigma t_1$ and $\sigma t_2$ as in the case of $\THH(\mathbb Z_p)/p$. For $p>2$, $\mu_2$ can be identified with $\sigma^2v_2$.
\end{exm}

\subsection{The associated graded}
We further filter the associated graded ring $j_{\zeta}^{\gr}$ by the $p$-adic filtration to ultimately reduce the computation to our understanding of $\THH(\FF_p)$. In running the spectral sequences to obtain the $\THH$ mod $(p,v_1)$, we find that they are close enough to the spectral sequences of $\ell_p^{B\ZZ}$, the fixed points of $\ell_p$ with the trivial $\ZZ$-action.

\begin{dfn}\label{defn:padicfil}
We define the $p$-adic filtration on $j_{\zeta}^{\gr}$ to be $j_{\zeta}^{\gr}\otimes_{\ZZ_p}\ZZ_p^{\fil}$. This is an $\mathbb E_{\infty}$-$\mathbb Z$-algebra object in the category of filtered graded spectra.

By taking the associated graded, we obtain $j_{\zeta}^{\gr}\otimes_{\mathbb Z_p}\mathbb F_p[v_0]$, which is an $\mathbb E_\infty$-$\mathbb Z$-algebra object in the category of bigraded spectra. We shall write \emph{hfp grading} for the grading on $j_{\zeta}^{\gr}$ if we need to distinguish it from the $p$-adic grading on $\mathbb F_p[v_0]$. For example, in $j_{\zeta}^{\gr}\otimes_{\ZZ_p}\FF_p[v_0]$, $v_1$ has hfp degree $2p-2$ and $p$-adic degree $0$, and $v_0$ has hfp degree $0$ and $p$-adic degree $1$.
\end{dfn}

\begin{lem}\label{lem:grdmodp}
	For $p>2$, there is an isomorphism of bigraded $\EE_{1}$-$\THH(\FF_p)$-algebras for $$\THH(j_{\zeta}^{\gr}\otimes_{\ZZ_p}\FF_p[v_0]) \cong \THH(\FF_p)\otimes_{\FF_p} \HoH(\FF_p[v_0,v_1]/\FF_p)\otimes_{\FF_p} \HoH(\FF_p^{B\ZZ}/\FF_p)$$
\end{lem}
\begin{proof}
	First note that $j_{\zeta}^{\gr}\otimes_{\ZZ_p}\FF_p[v_0] \cong j_{\zeta}^{\gr}/p\otimes_{\FF_p}\FF_p[v_0]$, which by \Cref{lem:modp} is equivalent to $\FF_p[v_1,v_0]\otimes_{\FF_p}\FF_p^{B\ZZ}$. Then, the statement follows from \Cref{lem:THHpoly}.
\end{proof}

We next study the behavior of fixed points by trivial $\ZZ$-actions on $\THH$. We use the spherical Witt vectors adjunction \cite[Proposition 2.2]{burklund2022chromatic} \cite[Section 5.2]{ellipticii} between perfect $\FF_p$-algebras and $p$-complete $\EE_{\infty}$-rings. For a perfect $\FF_p$-algebra $A$, $\WW(A)$ is an $\EE_{\infty}$-ring that is ($p$-completely) flat under $\SP_p$, and whose $\FF_p$ homology is $A$. The right adjoint is $\pi_0^{\flat}$ which is defined to be the inverse limit perfection of the $\FF_p$-algebra $\pi_0(R)/p$.

\begin{lem}\label{lem:spheretriv}
	There is an equivalence of $\EE_{\infty}$-$\SP_p^{B\ZZ}$-algebras $\THH(\SP_p^{B\ZZ}) \cong \SP_p^{B\ZZ}\otimes \WW(C^0(\ZZ_p;\FF_p))$. The restriction map $\SP_p^{B\ZZ} \to \SP_p^{Bp\ZZ}$ on $\pi_0^{\flat}$ is the map $C^0(\ZZ_p;\FF_p) \to C^0(p\ZZ_p;\FF_p)$ that restricts a function to $p\ZZ_p$.
\end{lem}
\begin{proof}
	There is a natural map $\SP^{B\ZZ}_p\otimes \WW(\pi_0^{\flat}(\THH(\SP_p^{B\ZZ}))) \to \THH(\SP_p^{B\ZZ})$, and so for the first claim it suffices to show that this is an equivalence and that $\pi_0^{\flat}(\THH(\SP_p^{B\ZZ})) \cong C^0(\ZZ_p;\FF_p)$. Both of these can be checked after base change to $\FF_p$. Note that $\THH(\SP_p^{B\ZZ})_p\otimes \FF_p\cong \HoH(\FF_p^{B\ZZ}/\FF_p)$.
	
	Since $\FF_p^{B\ZZ} = \colim_n \FF_p^{\text{B}\ZZ/p^n\ZZ}$ and $\text{B}\ZZ/p^n\ZZ$ is $p$-finite, we have, by \cite[Corollary 1.1.10]{DAGXIII}, $$\HoH(\FF_p^{\text{B}\ZZ/p^n\ZZ}/\mathbb F_p)\cong \FF_p^{\text{B}\ZZ/p^n\ZZ}\otimes_{\FF_p^{(\text{B}\ZZ/p^n\ZZ)^2}}\FF_p^{\text{B}\ZZ/p^n\ZZ} \cong \FF_p^{\text{B}\ZZ/p^n\ZZ\times_{(\text{B}\ZZ/p^n\ZZ)^2} \text{B}\ZZ/p^n\ZZ}.$$
	We have equivalences of spaces natural in $n$
	\[
	\text{B}\ZZ/p^n\ZZ \times_{(\text{B}\ZZ/p^n\ZZ)^2}\text{B}\ZZ/p^n\ZZ \cong \text{LB}\ZZ/p^n\ZZ = \text{B}\mathbb Z/p^n\ZZ \times \mathbb Z/p^n\ZZ.
	\]
	where L denotes the free loop space.
	Then, via the K\"unneth isomorphism and taking the colimit over $n$, we get \[
	\HoH(\FF_p^{B\ZZ}/\FF_p) \cong \colim_n \HoH(\FF_p^{B\ZZ/p^n\ZZ}/\FF_p) \cong\colim_n(\FF_p^{B\ZZ/p^n\ZZ}\otimes \FF_p^{\ZZ/p^n\ZZ}) \cong\FF_p^{B\ZZ}\otimes \colim_n\FF_p^{\ZZ/p^n\ZZ}.
	\]
	Since $\colim_n\FF_p^{\ZZ/p^n\ZZ}$ is $\ctf{\ZZ_p}$, we obtain the desired equivalence. 
	
	To see the claim about $\pi_0^{\flat}$, we note the natural map $\FF_p^{B\ZZ} \to \FF_p^{Bp\ZZ}$ is the colimit of  $\FF_p^{hB\ZZ/p^n\ZZ} \to \FF_p^{hB\ZZ/p^{n-1}\ZZ}$, where the map is given by the inclusion $\ZZ/p^{n-1}\ZZ \to \ZZ/p^n\ZZ$. At the level of the $\pi_0$, $LB\ZZ/p^{n-1}\ZZ \to LB\ZZ/p^n\ZZ$ is also the inclusion $\ZZ/p^{n-1}\ZZ \to \ZZ/p^n\ZZ$, so induces the restriction map at the level of $\ctf{-}$. Taking the colimit over $n$ gives the claim.
\end{proof}
\begin{rmk}
	\Cref{lem:spheretriv} can be interpreted as saying that the failure of $p$-adic $\THH$ to commute with taking $\ZZ$-homotopy fixed points in the universal case is measured by $\pi_0^{\flat}$. In particular, the map $\THH(\SP_p^{B\ZZ}) \xrightarrow{f} \THH(\SP_p)^{B\ZZ}$ on $\pi_0^{\flat}$ is the map $\ctf{\ZZ_p} \xrightarrow{\pi_0^{\flat}f} \FF_p$ evaluating at $0$, and the comparison map is base changed along $\WW(\pi_0^{\flat}f)$.
\end{rmk}

\begin{cor}\label{cor:trivfixed}
	Let $R$ be a $p$-complete $\EE_{\infty}$-ring. Then there is an equivalence of $p$-complete $\EE_{\infty}$-$R$-algebras $\THH(R^{B\ZZ})\cong \THH(R)^{B\ZZ}\otimes \WW(\ctf{\ZZ_p})$.
\end{cor}

\begin{ntn}
	We use $\zeta$ to refer to the class in $\pi_{-1}\SP^{B\ZZ}$ that is the generator of the cohomology of the circle.
\end{ntn}

We note that the map $\SP^{B\ZZ} \to j_{\zeta}$ sends $\zeta$ to the generator of $\pi_{-1}j_{\zeta}$. This explains the name, since the class $\zeta \in \pi_{-1}j_{\zeta}$ accounts for the difference between $j_{\zeta}$ and $j$.

Combining \Cref{cor:trivfixed} with \Cref{lem:grdmodp} and the HKR isomorphism, we get the following.
\begin{cor} For $p>2$, we have an isomorphism of rings
\[ 
	\pi_*\THH(j_{\zeta}^{\gr}\otimes_{\ZZ_p}\FF_p[v_0]) \cong \FF_p[\sigma^2p,v_0,v_1]\otimes\Lambda[dv_0,dv_1,\zeta]\otimes \ctf{\ZZ_p}
	\]
\end{cor}

\subsection{Spectral sequences}

Let us first run the spectral sequence for the $p$-adic filtration.

\begin{prop}\label{prop:THHjzetagr}
	For $p>2$, we have an isomorphism of rings
	\[
	\pi_*\THH(j_{\zeta}^{\gr})/p \cong \pi_*\THH(\ZZ_p)/p\otimes \FF_p[v_1]\otimes\Lambda[dv_1,\zeta]\otimes \ctf{\ZZ_p}.
	\]
\end{prop}

\begin{proof}
	As in Example \ref{exm:THHZ}, using \Cref{lem:spheretriv}, we find that the spectral sequence associated with $\THH(j_{\zeta}^{\gr}\otimes \ZZ_p^{\fil})/\tilde{v_0}$ has $E_1$-page isomorphic to $\pi_*\THH(j_{\zeta}^{\gr}\otimes_{\ZZ_p}\FF_p[v_0])/v_0 \cong \FF_p[\sigma^2p,v_1]\otimes\Lambda[dv_0,dv_1]\otimes \Lambda[\zeta] \otimes C^0(\ZZ_p;\FF_p)$ and converges to $\pi_\ast\THH(j_{\zeta}^{\gr})/p$.
	
	Because there is a map of filtered rings $j_{\zeta}^{\gr}\otimes_{\ZZ_p}\ZZ_p^{\fil} \to \THH(j_{\zeta}^{\gr}\otimes_{\ZZ_p}\ZZ_p^{\fil})$, we see that the classes $v_1, \zeta$ are permanent cycles. The class $dv_1$ is a permanent cycle since it detects the suspension $dv_1$ of $v_1\in\pi_\ast j_{\zeta}^{\gr}/p$. The elements of $C^0(\ZZ_p;\FF_p)$ are permanent cycles since there are no elements of negative topological degree and positive filtration.
	
	From the map of filtered rings
	\[
	\THH(\mathbb Z_p^{\fil})\to \THH(j_{\zeta}^{\gr}\otimes_{\ZZ_p}\ZZ_p^{\fil}),
	\]
	there is a $d_1$-differential $\sigma^2p\mapsto \sigma v_0$ by Example \ref{exm:THHZ}, and $(\sigma^2p)^p$ and $(\sigma^2p)^{p-1}dv_0$ are permanent cycles detecting images of classes in $\THH(\ZZ_p)$. It follows that after the $d_1$-differential, the $E_2$-page is $\FF_p[(\sigma^2p)^p,v_1]\otimes\Lambda[(\sigma^2p)^{p-1}dv_0,dv_1,\zeta]\otimes C^0(\ZZ_p;\FF_p)$, so the spectral sequence collapes at the $E_2$-page. There are no multiplicative extensions since every class comes from either $j_{\zeta}^{\gr}$, $\THH(\mathbb Z_p)$, or $\THH(\mathbb S_p^{B\mathbb Z})$.
\end{proof}

Our next goal is to compute mod $(p,v_1)$ the spectral sequence $\THH(j^{\gr}_{\zeta}) \implies \THH(j_{\zeta})$. 
Before doing so, we run the analogous spectral sequence for computing $\THH(\ell_p)/(p,v_1)$, as a warm up. We consider the $\EE_{\infty}$-ring $\ZZ_{\zeta}:=\ZZ_p^{B\ZZ}$ with the trivial filtration.

\begin{thm}\label{thm:jzetathh}
	For $p>2$, $\pi_*(\THH(j_{\zeta}))/(p,v_1) \cong \FF_p[\sigma^2v_2]\otimes\Lambda[\lambda_1,\lambda_2,\zeta]\otimes C^0(\ZZ_p;\FF_p)$ with $|\lambda_i| = 2p^i-1$ and $|\sigma^2v_2| = 2p^2$.
\end{thm}

\begin{proof}
	\begin{figure}[t]
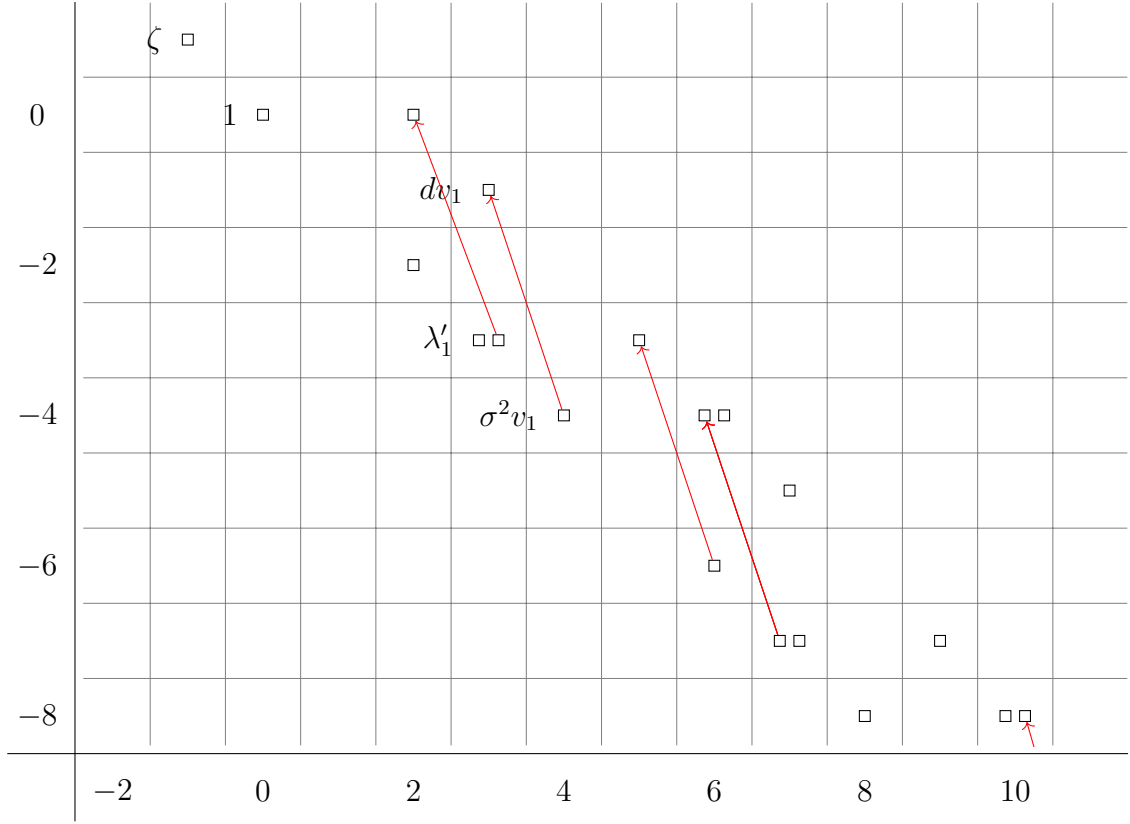

		\centering
		\scalebox{1.00}{
			\printpage[ name = juzeta2, page = 2 ]
		}
		\caption{
			Above is the spectral sequence associated with the filtered ring $\THH(ju^{\fil})/(2,v_1)$. This spectral sequence is the $p=2$ version of the spectral sequence in \Cref{thm:jzetathh} (see \Cref{thm:juk}), and only has $d_2$ differentials. Each square represents a copy of $C^0(\ZZ_2;\FF_2)$.
		}
		\label{fig:thhju}    
	\end{figure}
	As in Example \ref{exm:THHell}, we consider the spectral sequence associated with the filtered spectrum $\THH(j_{\zeta}^{\fil})/(p,\tilde{v_1})$. The analogous spectral sequence in the case $p=2$ is displayed in \Cref{fig:thhju} above. The underlying spectrum is $\THH(j_{\zeta})/(p,v_1)$ and the associated graded spectrum is $\THH(j_{\zeta}^{\gr})/(p,v_1)$. By \Cref{prop:THHjzetagr}, the $E_1$-page is isomorphic to $\mathbb F_p[\sigma^2 v_1]\otimes\Lambda[\lambda_1, dv_1,\zeta]\otimes C^0(\mathbb Z_p;\mathbb F_p)$.

 The classes in $C^0(\mathbb Z_p;\mathbb F_p)$ are permanent cycles by the Leibniz rule, since they are all their own $p^{th}$-power. The class $\zeta \in H^1(S^1;\FF_p)$ is a permanent cycle because it detects a class in the image of $j_{\zeta}\to \THH(j_{\zeta})$.

	By \Cref{lem:diffonsigma}, there is a differential $d_{2p-2}(\sigma^2 v_1)\doteq dv_1$, and the Leibniz rule determines the differentials on powers of $\sigma^2v_1$. 
 
        Similarly, by Lemma \ref{lem:diffonsigma}, there must be a $d_{2p-2}$ differential $\lambda_1\doteq \sigma^2 \alpha_1\dotmapsto d\alpha_1$ in the spectral sequence $\THH(j_{\zeta}^{\gr})\implies \THH(j_\zeta)$ mod $p$. By Lemma \ref{lem:dLeibniz}, we have
	\[
		d\alpha_1 = d(v_1\zeta) = v_1d\zeta - \zeta dv_1,
	\]
	so that we have the differential $d_{2p-2}(\lambda_1) \doteq\zeta dv_1$ mod $(p,v_1)$. By using the previous paragraph and replacing $\lambda_1$ with
 \[
    \lambda_1'=\lambda_1 - \epsilon \zeta\mu
 \]
 for some $\epsilon\in\mathbb F_p^\times$, we may assume that $d_{2p-2}(\lambda_1')=0$.
 	
	This completely determines the spectral sequence up to the $E_{2p-2}$-page, and we learn the $E_{2p-1}$-page is isomorphic to $\mathbb F_p[(\sigma^2v_1)^p]\otimes\Lambda[\lambda_1',(\sigma^2v_1)^{p-1}dv_1,\zeta]\otimes C^0(\ZZ_p;\mathbb F_p)$. There are no more differentials since there is no class outside filtration degree $0$ and $2p-2$.
 There are no multiplicative extension problems since the multiplicative generators in nonzero degree are free generators as a graded ring. 
	
	Finally, let us show that the polynomial generator $\mu_2$ is the class $\sigma^2 v_2$. Let us consider the map $j_{\zeta}^{\fil} \to \ZZ_{\zeta}$ induced by applying $(\tau_{\geq*}(-))^{B\ZZ}$ to the $\ZZ$-equivariant truncation map $\ell_p \to \ZZ_p$. This induces a map of spectral sequences for $\THH$. Since $\ZZ_\zeta$ has the trivial filtration, its $\THH$ does too, so has no differentials in its associated spectral sequence. By \Cref{cor:trivfixed} and \Cref{exm:THHZ}, $\THH(\ZZ_{\zeta})_p\cong \THH(\ZZ_p)^{B\ZZ}\otimes \WW(\ctf{\ZZ_p})$, so $$\pi_*\THH(\ZZ_{\zeta})/(p,v_1) \cong \FF_p[\sigma^2v_1]\otimes \Lambda[\lambda_1,\zeta]\otimes \ctc{\ZZ_p}.$$
		
	$v_2 \in \pi_{2p^2-2}\SP/(p,v_1)$ has a canonical nullhomotopy in $j_{\zeta}/(p,v_1) \cong \FF_p^{B\ZZ}$ and $\ZZ_{\zeta}/(p,v_1) \cong \FF_p[\sigma v_1]^{B\ZZ}$, so there is a canonical element $\sigma^2v_2$ in $\pi_{2p^2}\THH(j_{\zeta})/(p,v_1)$ and $\pi_{2p^2}\THH(\ZZ_{\zeta})/(p,v_1)$, which we claim is detected in the spectral sequence for $\THH(j_{\zeta}^{\fil})$ by $(\sigma^2v_1)^p$. To see this, it suffices to show this in $\THH(\ZZ_{\zeta})$ because the map is injective in degree $2p^2-2$. But now it is the image of $\sigma^2v_2$ from the map $\ell_p^{\fil} \to \ZZ_{\zeta}$, and in $\ell_p^{\fil}$, which we know by \Cref{exm:THHell} is detected by $(\sigma^2v_1)^p$.
\end{proof}

\begin{rmk}\label{rmk:lambda1}
In the proof of the previous theorem, a reader might wonder why $\lambda_1$ supports a differential while $\sigma^2\alpha_1$ is still well-defined in $\THH(j_{\zeta})$. This can be explained by the fact that $\sigma^2\alpha_1$ is not well-defined in $\THH(\mathbb Z_{\zeta})/(p,v_1)$ since
	\[
		\pi_{2p-3}(\fib(\mathbb S\to\mathbb Z_{\zeta})/(p,v_1)) \to \pi_{2p-3}(\mathbb S/(p,v_1))
	\]
	is not injective. The class $\sigma^2\alpha_1$ is well-defined in $\THH(\mathbb Z)/(p,v_1)$ and $\THH(j_{\zeta})/(p,v_1)$, but their images in $\THH(\mathbb Z_{\zeta})/(p,v_1)$ are different. The class $\lambda_1$ in the $E_1$-page represents the former and $\lambda_1'$ represents the latter.
\end{rmk}

\begin{rmk}\label{rmk:THHellhZ}
	We can carry out the same computation for $\THH(\ell_p)^{h\mathbb Z}/(p,v_1)$ using the same filtrations $\ell_p^{\fil}$ and $\ell_p^{\gr}\otimes \ZZ_p^{\fil}$. Then, we obtain an isomorphism of rings
	\[
		\pi_\ast\THH(\ell_p)^{h\mathbb Z}/(p,v_1) \simeq \mathbb F_p[\sigma^2 v_2] \otimes\Lambda[\lambda_1,\lambda_2,\zeta].
	\]
	Furthermore, by keeping track of the map
	\[
		\THH(j_{\zeta})/(p,v_1)\to\THH(\ell_p)^{h\mathbb Z}/(p,v_1)
	\]
	at every stage, we see that on homotopy groups, this map is the base-change along
	\[
		C^0(\mathbb Z_p;\mathbb F_p)\to \mathbb F_p
	\]
	that evaluates a function at $0\in\mathbb Z_p$. The key point is in the proof of \Cref{lem:spheretriv}, where one uses the fact that the natural map $B\ZZ/p^n \to LB\ZZ/p^n$ coming from constant loops is the inclusion of the component $0$ for each $n\geq 0$.
\end{rmk}

\subsection{The prime 2}
We next turn to the prime $2$. We first need to run the analogous analysis as in \Cref{exm:THHell} for $\ko_2$. We consider $\ko_2^{\gr}/2[v_0]$ as the bigraded ring given as the associated graded of $\ko_2^{\gr}\otimes_{\ZZ_2}\ZZ_2^{\fil}$. To understand this, we need the following lemma.

\begin{lem}\label{lem:thhkogrmod2}
	There is an isomorphism of bigraded rings $$\pi_*\THH(\ko_2^{\gr}/2[v_0])/\eta \cong \FF_2[v_0,v_1,\sigma^22,d\eta]/((d\eta)^2+v_1d\eta) \otimes \Lambda[dv_0,dwv_1]$$
\end{lem}

\begin{proof}
	The associated graded of $\ko_2^{\gr}/2[v_0]$ with respect to the Posnikov filtration is $\FF_2[v_0,v_1,\eta]$. 
	
	By symmetric monoidality of $\THH$, we have an equivalence
	\[
	\THH(\FF_2[v_0,v_1,\eta])\cong \THH(\FF_2[v_0,v_1])\otimes_{\THH(\FF_2)}\THH(\FF_2[\eta])
	\]
	Since the argument of \Cref{lem:grdmodp} works at the prime $2$, we learn that the first tensor factor has homotopy ring $\FF_2[\sigma^2p,v_0,v_1]\otimes \Lambda[dv_0,dv_1]$.
	
	For the second tensor factor, we note that $\THH(\FF_2[\eta])\otimes_{\THH(\FF_2)}\FF_2 \cong \HoH(\FF_2[\eta]/\FF_2)$, whose homotopy ring is $\FF_2[\eta]\otimes \Lambda[d\eta]$. Since the map $\THH(\FF_2) \to \FF_2$ is the cofiber of $\sigma^2p$, we can run a $\sigma^2p$-Bockstein spectral sequence to recover $\THH(\FF_2[\eta])$. In the spectral sequence, $\eta,d\eta$ are permanent cycles since they are in the image of the unit map and the map $d$. We also see that there are no multiplicative extensions mod $\eta$ for degree reasons, i.e. we have
	\[
		\pi_\ast\THH(\mathbb F_2[\eta])/\eta = \Lambda(d\eta)\otimes_{\mathbb F_2}\mathbb F_2[\sigma^22].
	\]
	
In the spectral sequence computing $\THH(\ko_2^{\gr}/2[v_0])$ from this, everything is a permanent cycle since all classes are generated either from the image of the unit map, the map from $\THH(\FF_2)$, or the map $d$.
	
	Now we turn to the multiplicative extensions, which we compute by mapping to the $\sigma^22$-completion of $\THH(\FF_2^{BC_2}[v_0,v_1])$. As before, we can compute this via the $\sigma^22$-Bockstein spectral sequence whose $E_{1}$-page is $\HoH(\FF_2^{BC_2}[v_0,v_1]/\FF_2)[\sigma^22]$.
	
	We have an isomorphism $\HoH(\FF_2^{BC_2}[v_0,v_1]/\FF_2) \cong \HoH(\FF_2^{BC_2}/\FF_2)\otimes_{\FF_2}\HoH(\FF_2[v_0,v_1])$.
	Moreover, $\HoH_*(\FF_2[v_0,v_1]) \cong \FF_2[v_0,v_1]\otimes\Lambda[dv_0,dv_1]$, and $\HoH(\FF_2^{BC_2})$ is $\FF_2^{BC_2}\times \FF_2^{BC_2}$, since the free loop space of $BC_2$ is $BC_2\times C_2$. If $h$ is the generator of $\pi_{-1}\FF_2^{BC_2}$, then a nontrivial idempotent in $\pi_0\HoH(\FF_2^{BC_2})$ is given by $dh$. By the Leibniz rule (\Cref{lem:dLeibniz}), $d\eta = v_1dh+hdv_1$, so $(d\eta)^2 = v_1^2dh = v_1d\eta+\eta dv_1$. This this relation happens in $\THH(\ko_2)/\sigma^22$, but for degree reasons, this forces it to happen in $\THH(\ko_2)/\eta$ as well. 
	
	To see that the classes $dv_0$ and $dv_1$ square to $0$, we note that this is true in $\HoH(\FF_2[v_0,v_1]/\FF_2)$, and that we have a map
	$$\THH(\FF_2)\otimes_{\FF_2}\HoH(\FF_2[v_0,v_1]/\FF_2) \cong \THH(\FF_2[v_0,v_1]) \to \THH(\FF_2[v_0,v_1]^{BC_2})$$ using the isomorphism of \Cref{lem:THHpoly}.
\end{proof}

\begin{lem}\label{lem:thhkogr}
	There is an isomorphism of graded rings $$\pi_*(\THH(\ko_2^{\gr})/(2,\eta)) \cong \FF_2[v_1,\sigma^2v_1,d\eta]/((d\eta)^2+v_1d\eta)\otimes \Lambda[\sigma^2\eta,dv_1]$$
\end{lem}
\begin{proof}
	We now understand the spectral sequence computing $\pi_*(\THH(\ko_2^{\gr})/(2,\eta))$ by running the $2$-adic filtration spectral sequence on $\THH(\ko_2\otimes_{\ZZ_2}\ZZ_2^{\fil})/(v_0,\eta)$. By \Cref{exm:THHZ}, there is a differential from $\sigma^2v_0$ to $dv_0$, $\sigma^2\eta$ is a class squaring to zero detected by $\sigma^2v_0dv_0$, and $\sigma^2v_1$\footnote{The element $v_1\in \pi_2\SP/2$ exists, even though it does not extend to a self map.} detects $(\sigma^2v_0)^2$. The remaining classes are either in the image of the unit map or the image of $d$, so are permanent cycles. The relation $(d\eta)^2+v_1d\eta=0$ occurs because it does on associated graded, and because there are no classes in topological degree $4$ and positive $p$-adic filtration. The class $dv_1$ squares to zero since there are no classes of weight $-2$, topological degree $6$, and positive $p$-adic filtration.
\end{proof}

We now compute $\THH(\ko_2)/(2,\eta,v_1)$, which was also computed in \cite[Theorem 8.14]{angeltveit2005hopf}.

\begin{figure}[t]
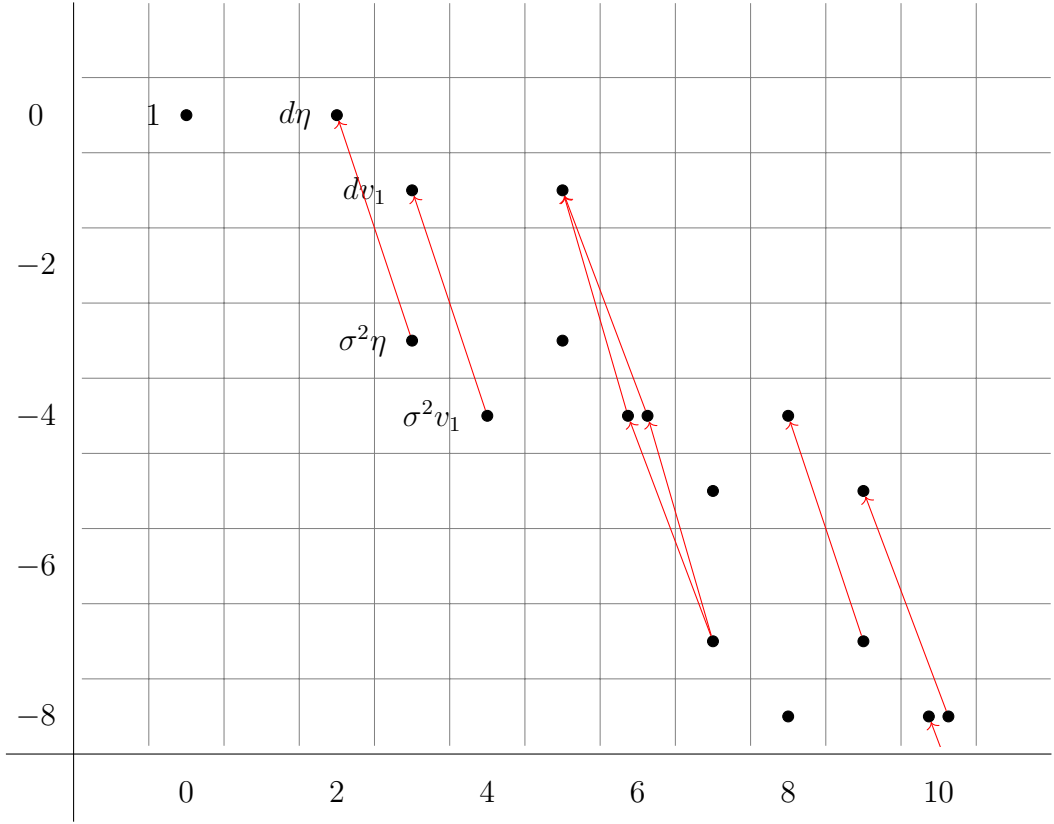

	\centering
	\scalebox{1.00}{
		\printpage[ name = thhko, page = 2 ]
	}
	\caption{
		Above is the spectral sequence associated with the filtered ring $\THH(\ko_2^{\fil})/(2,\eta,v_1)$.
	}
	\label{fig:thhko}    
\end{figure}

\begin{exm}\label{exm:thhko}
	We now can run the spectral sequence $$\THH(\ko_2^{\gr})/(2,\eta,v_1) \implies \THH(\ko_2)/(2,\eta,v_1)\text{,}$$ which is a spectral sequence associated with a filtered $\EE_{\infty}$-ring since $\FF_2 \cong \ko_2^{\fil}/(2,\eta,v_1)$, where $\eta$ and $v_1$ are taken in filtration $2$. This spectral sequence is displayed in \Cref{fig:thhko}. The first page of this spectral sequence by \Cref{lem:thhkogr} is $\FF_2[\sigma^2v_1]\otimes \Lambda[dv_1,d\eta,\sigma^2\eta]$. It follows as in \Cref{exm:THHell} that there are differentials from $\sigma^2\eta$ to $d\eta$ and $\sigma^2v_1$ to $dv_1$. What remains after these differentials are $\FF_2[(\sigma^2v_1)^2]\otimes \Lambda[\sigma^2v_1dv_1, \sigma^2\eta d\eta]$. For degree reasons, there can be no further differentials. the classes in odd degree square to $0$ because there are no classes in degrees $2$ or $6$ mod $8$.
\end{exm}

We now run the analogous analysis to compute $\THH(j_{\zeta})/(2,\eta,v_1)$.

\begin{lem}
	There is an isomorphism of graded rings $$\pi_*\THH(j_{\zeta}^{\gr})/(2,\eta,v_1) \cong \pi_*\THH(\ko_2^{\gr})/(2,\eta,v_1)\otimes \pi_*(\HoH(\FF_2^{B\ZZ}/\FF_2))$$
\end{lem}

\begin{proof}
	Since $\ko_2^{\gr}/2\otimes_{\FF_2}\FF_2^{B\ZZ} \cong j_{\zeta}^{\gr}/2$, we learn from \Cref{lem:thhkogrmod2} that $$\pi_*(\THH(j_{\zeta}^{\gr}/2[v_0])/(v_0,\eta,v_1)\cong \FF_2[\sigma^2v_0]\otimes \Lambda[d\eta,dv_0,dv_1]\otimes \HoH_*(\FF_2^{B\ZZ}/\FF_2)$$
	where $\HoH(\FF_2^{B\ZZ}/\FF_2)$ is computed via \Cref{cor:trivfixed} as $\FF_2^{B\ZZ}\otimes C^0(\mathbb Z_2;\mathbb F_2)$. 
	
	Exactly as in \Cref{lem:thhkogr}, in the spectral sequence for the $2$-adic filtration, there is a differential from $\sigma^2v_0$ to $dv_0$, $\sigma^2\eta$ is a class squaring to zero detected by $\sigma^2v_0dv_0$, and $\sigma^2v_1$ is a class detecting $(\sigma^2v_0)^2$. The rest of the classes are permanent cycles because they are either in the unit map, come from $d$, or are permanent cycles by the Leibniz rule. 
\end{proof}

\begin{thm}\label{thm:thhjzeta2}
	There is an isomorphism of rings for $p=2$ $$\pi_*\THH(j_{\zeta})/(2,\eta,v_1) \cong \FF_2[\mu]\otimes \Lambda[\lambda_2,x,\zeta]\otimes C^0(\mathbb Z_2;\mathbb F_2)$$ where $|x| = 5$, $|\lambda_2| = 7, |\mu| = 8$.
\end{thm}

\begin{proof}
	We run the spectral sequence $\THH(\ko_2^{\gr})/(2,\eta,v_1) \implies \THH(\ko_2)/(2,\eta,v_1)$.
	As in \Cref{exm:thhko}, there are differentials from $\sigma^2\eta$ to $d\eta$ and $\sigma^2v_1$ to $dv_1$. For degree reasons, $(\sigma^2v_1)^2$ is a permanent cycle, as are $\sigma^2\eta d\eta$, $\sigma^2v_1dv_1$, and $\zeta$. All classes in $C^0(\mathbb Z_2;\mathbb F_2)$ are permanent cycles by the Leibniz rule. If we let $\lambda_2$ and $x$ denote classes detecting $\sigma^2v_1dv_1$ and $\sigma^2\eta d\eta$ respectively, then $\lambda_2^2=0$ and $x^2=0$ for degree reasons. 
\end{proof}

\section{The $\THH$ of $j$}\label{sec:j}
We now consider $\THH(j)/(p,v_1)$ for $p>2$. We first compute the Hochschild homology of the $\FF_p$-algebra $j^{\gr}/p$, which is isomorphic to $\tau_{\geq0}(\FF_p[v_1]^{B\ZZ})$ by \Cref{lem:modp}.

\begin{prop}\label{prop:hhjgr/p}
	Let $p>2$. $\HoH_*((j^{\gr}/p)/\FF_p) \cong \HoH_*(\tau_{\geq0}(\FF_p[v_1]^{B\ZZ})/\FF_p)$ is isomorphic as a ring to 
	\[
		\Lambda[dv_1,\alpha_1]\otimes \mathbb F_p[v_1,x_0,x_1,\dots]/(x_i^p = v_1^{p^{i+1}-p^i}x_i + v_1^{p^{i+1}-p^i - 1} \alpha_1 (\prod_{j=0}^{i-1}x_j^{p-1})dv_1;\;i\geq0)
	\]
	where $|x_i|=p^i(2p-2)$, and $x_i$ is in grading $p^i(2p-2)$.
\end{prop}
\begin{proof}
	Define a graded ring $R = \tau_{\geq 0}\mathbb Z_p[v_1]^{B\mathbb Z}$ so that $R/p\simeq j^{\gr}/p$. We shall show that $\pi_\ast\HoH(R/\mathbb Z_p)$ is the $\mathbb Z_p$-algebra generated by $v_1,dv_1,\alpha$, and a set of generators $x_0,x_1,\dots$ with $|x_i| = p^i(2p-2)$ having relations
	\[
		x_i^p = px_{i+1} + v_1^{p^{i+1}-p^i}x_i + v_1^{p^{i+1}-p^i - 1} \alpha (\prod_{j=0}^{i-1}x_j^{p-1})dv_1.
	\]
	Then, the statement follows by the base-change $\mathbb Z_p\to\mathbb F_p$.
	
	Let $R_{\zeta} = \mathbb Z_p[v_1]^{B\mathbb Z}$ and let $\eta:R\to R_{\zeta}$ denote the connective cover map. To compute the Hochschild homology, we shall show that the map $\eta_\ast:\pi_\ast\HoH(R)\to\pi_\ast\HoH(R_{\zeta})$ is injective and describe the image. Note that $\pi_\ast R_{\zeta} = \mathbb Z_p[v_1,\zeta]$ and $\pi_\ast R = \mathbb Z_p[v_1,\alpha]$ where $\eta_\ast(\alpha)= v_1\zeta$.
	
	Let us consider the K\"unneth spectral sequence
	\begin{equation}\label{eq:KssR}
		E_2(\HoH(R)) = \Tor^{\pi_\ast(R\otimes_{\mathbb Z} R)}(\pi_\ast R,\pi_\ast R) \implies \pi_\ast \HoH(R).
	\end{equation}
	Since $\pi_\ast R =\mathbb Z_p[v_1]\otimes\Lambda[\alpha]$, the $E_2$-page can be computed as
	\[
		E_2(\HoH(R)) = \mathbb Z_p[v_1]\otimes\Lambda[dv_1,\alpha]\otimes\Gamma[d\alpha].
	\]
	Similarly, there is a spectral sequence
	\begin{equation}\label{eq:KssRzeta}
		E_2(\HoH(R_{\zeta})) = \mathbb Z_p[v_1]\otimes\Lambda[dv_1,\zeta]\otimes\Gamma[d\zeta]\implies \pi_\ast \HoH(R_{\zeta})
	\end{equation}
	up to $p$-completion.
	
	We claim that $E_2(\HoH(R))\to E_2(\HoH(R_{\zeta}))$ is injective. By Lemma \ref{lem:dLeibniz}, we have
	\[
		d\alpha \mapsto -\zeta dv_1 + v_1d\zeta.
	\]
	To prove the injectivity, it is enough to prove it after taking the associated graded group with respect to the $(dv_1)$-adic filtration. Then, we may assume that $d\alpha$ maps to $v_1d\zeta$, and since $E_2(\HoH(R_{\zeta}))$ is torsion-free, the divided power $\gamma_n(d\alpha)$ maps to $v_1^n\gamma_n(d\zeta)$. Therefore, we have the desired injectivity. Note also that the map is injective mod $p$.
	
	The spectral sequence \eqref{eq:KssRzeta} degenerates at the $E_2$-page using the symmetric monoidality of $\HoH$, \Cref{cor:trivfixed}, and \Cref{lem:hkrdx}. We then see that \eqref{eq:KssR} also degenerates at the $E_2$-page and that $\eta_\ast:\pi_\ast\HoH(R)\to\pi_\ast\HoH(R_{\zeta})$ is injective, even after mod $p$.
	
	Let us describe the K\"unneth filtration on
	\[
		\pi_\ast\HoH(R_\zeta) = \mathbb Z_p[v_1]\otimes\Lambda[dv_1,\zeta]\otimes W(C^0(\mathbb Z_p;\mathbb F_p))
	\]
	in more detail. Here, the ring
	\[
	W(C^0(\mathbb Z_p;\mathbb F_p))=\lim_{k}C^0(\mathbb Z_p;\mathbb Z_p/p^k)
	\] 
	is the ring of all continuous functions $\mathbb Z_p\to\mathbb Z_p$. It can also be described, up to completion, as the algebra generated by $y_0,y_1,\dots$ with relations
	\[
		y_i^p = py_{i+1}+y_i.
	\]
	Here, the element $y_0$ is the identity function $\mathbb Z_p\to\mathbb Z_p$ and the $y_i$'s for $i>0$ can be defined with the above formula since $y^p\equiv y\pmod p$ for any $y\in W(C^0(\mathbb Z_p;\mathbb F_p))$. In $\pi_\ast\HoH(R_\zeta)$, the element $y_0$ equals $d\zeta$, and the $y_i$'s represent the $p^i$-th divided power of $d\zeta$ in the K\"unneth spectral sequence \eqref{eq:KssRzeta}.
	
	To determine $\pi_\ast\HoH(R)$, we need to find the classes $x_i$'s representing the divided powers $\gamma_{p^i}(d\alpha)\in E_2(\HoH(R))$ up to a $p$-adic unit. The first divided power $d\alpha\in E_2(\HoH(R))$ has a canonical lift $x_0:=d\alpha\in\pi_\ast\HoH(R)$ and its image under $\eta_\ast$ is $v_1y_0 - \zeta dv_1$. Inductively, suppose that we have chosen $x_0,\dots,x_i$ in a way that the image of $x_j$ is
	\[
		\eta_\ast(x_j)= v_1^{p^j}y_j - v_1^{p^j-1}(\prod_{k=0}^{j-1}y_k^{p-1})\zeta dv_1
	\]
	for $0\leq j\leq i$. Let $x_{i+1}$ be any class representing $\gamma_{p^{i+1}}(d\alpha)$. Then, after scaling by a unit, we must have
	\[
		x_i^p = px_{i+1}+c
	\]
	for some class $c\in \pi_\ast\HoH(R)$ with K\"unneth filtration $<p^{i+1}$. Applying $\eta_\ast$, we have
	\[
		\eta_\ast(c)\equiv \eta_\ast(x_i)^p \equiv v_1^{p^{i+1}}y_i^p \equiv v_1^{p^{i+1}}y_i \pmod p. 
	\]
	Let $d\in\pi_\ast\HoH(R)$ be the class $v_1^{p^{i+1}-p^i-1}(v_1x_i + \alpha(\prod_{k=0}^{i-1}x_k^{p-1})dv_1)$, having K\"unneth filtration $p^i$. Then, we can compute that $\eta_\ast(d) = v_1^{p^{i+1}}y_i$ so that $\eta_\ast(c) \equiv \eta_\ast(d) \pmod p$.
	Since $\eta_\ast$ is injective mod $p$, we have $c\equiv d\pmod p$, so by replacing $x_{i+1}$ with $x_{i+1} - (c - d)/p$, we can assume that $c=d$. Then, we have
	\begin{align*}
		\eta_\ast(x_{i+1}) &= p^{-1}\eta_\ast(x_i^p - c) \\&= p^{-1}\left( v_1^{p^{i+1}}y_i -pv_1^{p^{i+1}-1}(y_i\cdots y_0)^{p-1}\zeta dv_1 - v_1^{p^{i+1}}y_i\right)\\&=v_1^{p^{i+1}}y_{i+1} - v_1^{p^{i+1}-1}(y_i\cdots y_0)^{p-1}\zeta dv_1.
	\end{align*}
	The desired ring structure of $\pi_\ast\HoH(R)$ can now be read off from the ring structure on $\pi_*\HoH(R_{\zeta})$.
\end{proof}

\begin{lem}\label{lem:jgrdmodp}
	There is an isomorphism of bigraded $\EE_{1}$-$\THH(\FF_p)$-algebras for $p>2$ $$\THH(j^{\gr}\otimes_{\ZZ_p}\FF_p[v_0]) \cong \THH(\FF_p)\otimes_{\FF_p} \HoH(\FF_p[v_0]/\FF_p)\otimes_{\FF_p} \HoH(\tau_{\geq0}\FF_p[v_1]^{B\ZZ}/\FF_p)$$
\end{lem}
\begin{proof}
	We run the strategy of \Cref{lem:grdmodp} with appropriate modifications.  
	First, we have the isomorphism $j^{\gr}\otimes_{\ZZ_p}\FF_p[v_0] \cong j^{\gr}/p\otimes_{\FF_p}\FF_p[v_0]$, which by \Cref{lem:modp} is equivalent to $\tau_{\geq0}\FF_p[v_1,v_0]\otimes_{\FF_p}\FF_p^{B\ZZ}$. As an $\EE_2$-ring, we claim this is equivalent to the tensor product of $\FF_p\otimes \SP[v_0]$ with the pullback of the cospan \begin{center}
		\begin{tikzcd}
			{} & \SP[v_1]\otimes \SP^{B\ZZ} \ar[d]\\
			\SP\ar[r] & \SP^{B\ZZ}
		\end{tikzcd}
	\end{center}
	where the vertical map is the augmentation sending $v_1$ to $0$.
	
	This isomorphism is a consequence of the isomorphism of \Cref{lem:grdmodp} and the pullback square
	
	\begin{center}
		\begin{tikzcd}
			j^{\gr}/p \ar[r]\ar[d] & j_{\zeta}^{\gr}/p \ar[d]\\
			\ar[r] \FF_p& \FF_p^{B\ZZ}
		\end{tikzcd}
	\end{center}
	Given this equivalence, we conclude by arguing exactly as in \Cref{lem:grdmodp}.
\end{proof}

\begin{prop}\label{prop:thhjgr}
	Let $p>2$. Then
	\[
	\pi_*\THH(j^{\gr})/p \cong \pi_*\THH(\ZZ_p)/p\otimes \pi_*\HoH(\tau_{\geq0}\FF_p[v_1]^{B\ZZ}/\FF_p)
	\]
\end{prop}
\begin{proof}
	We follow the strategy in \Cref{prop:THHjzetagr}, running the spectral sequence corresponding to the $p$-adic filtration
	$$\pi_*\THH(j^{\gr}/p[v_0])/p \implies \pi_*\THH(j^{\gr})/p. $$
	The $E_1$-page is understood via \Cref{lem:jgrdmodp} to be 
	$$\FF_p[\sigma^2p,v_0]\otimes \Lambda[dv_0]\otimes \pi_*\HoH(\tau_{\geq0}\FF_p[v_1]^{B\ZZ}/\FF_p)$$
	where the last tensor factor is described in \Cref{prop:hhjgr/p}. There is a differential $d_1\sigma^2p = dv_0$, coming from the map from $\ZZ_p^{\fil} \to \ZZ_p^{\fil}\otimes j^{\fil}$ and \Cref{exm:THHZ}.
	
	We need to show that the remaining classes are permanent cycles. The classes $v_1,\alpha$ are permanent cycles because they are in the image of the unit map, and $dv_1$ is a permanent cycle because it is in the image of the map $\sigma^2$. The classes $x_i$ are permanent cycles for degree reasons, as everything of positive $p$-adic filtration is in nonnegative degree, and the differentials respect the hfp grading. One also sees for degree reasons and the map from $\THH(\ZZ_p)/p$ that there are no multiplicative extension problems.
\end{proof}

We now run the spectral sequence $\THH(j^{\gr})/(p,v_1) \implies \THH(j)/(p,v_1)$ associated with the filtered spectrum $\THH(j^{\fil})/(p,\widetilde{v_1})$ where $\widetilde{v_1}\in\pi_\ast j/p$ is the class of filtration $2p-2$. The following lemma guarantees the multiplicativity of the spectral sequences.

\begin{lem}\label{lem:todamult}
	$j^{\fil}/(p,\widetilde{v_1})$ admits a homotopy commutative $\AA_{p-1}$-multiplication for $p>2$, and in particular is homotopy associative for $p>3$.
\end{lem}

\begin{proof}
	By \cite[Example 3.3]{angeltveit2008topological}, it follows that $\SP/p$ is an $\AA_{p-1}$-algebra, and it is easy to see that there is no obstruction to its multiplication being homotopy commutative for $p>2$. We conclude by observing that $j^{\fil}/(p,\widetilde{v_1}) \cong \tau_{\leq 2p-3}j^{\fil}\otimes \SP/p$.
	
	Note that by loc. cit., the multiplication is not $\AA_p$, the obstruction being $\alpha_1$.
\end{proof}

\begin{figure}[t]
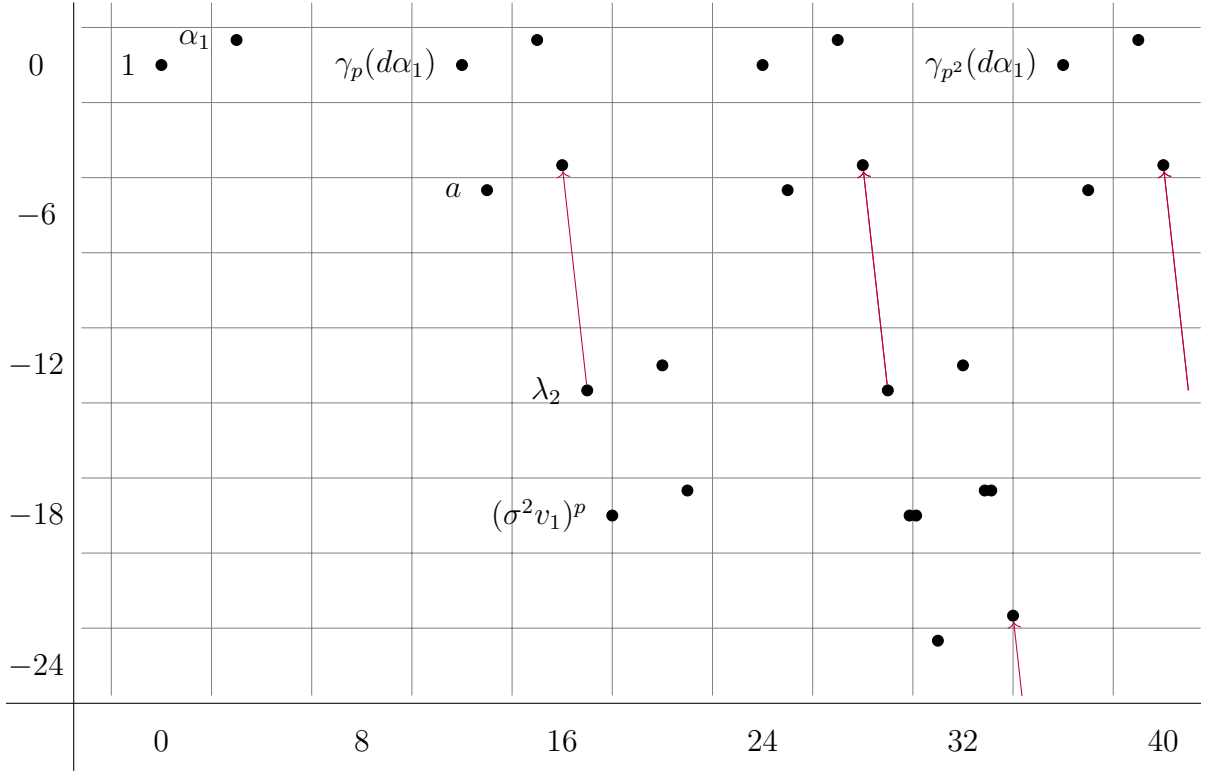

	\centering
	\scalebox{1.00}{
		\printpage[ name = thhj, page = 2 ]
	}
	\caption{
		Above is the $E_8$-page of the spectral sequence associated with the filtered ring $\THH(j^{\fil})/(3,v_1)$. The spectral sequence collapses at the $E_9$-page.
	}
	\label{fig:thhj}    
\end{figure}

\begin{thm}\label{thm:thhj}
	For $p>3$, $\pi_*\THH(j)/(p,v_1)$ is the homology of the CDGA
	\[
	\mathbb F_p[\mu_2]\otimes\Lambda[\alpha_1,\lambda_2,a]\otimes\Gamma[b],\quad d(\lambda_2)=a\alpha_1
	\]
	\[ |b| = 2p^2-2p , \; |a| = 2p^2-2p-1,\; |\lambda_2| = 2p^2-1,\; |\mu_2| = 2p^2\]
	and for $p=3$, the above result is true after taking an associated graded ring.
\end{thm}
\begin{proof}
	The $E_1$-page of the spectral sequence
	\[
		E_1 = \pi_\ast\THH(j^{\gr})/(p,v_1)\implies \pi_\ast\THH(j)/(p,v_1)
	\]
	is isomorphic to
	\[
		\mathbb F_p[\mu_1]\otimes\Lambda[\sigma^2\alpha_1,dv_1,\alpha_1]\otimes\Gamma[d\alpha_1].
	\]
	by \Cref{prop:thhjgr}. By \Cref{lem:diffonsigma}, there are $d_{2p-2}$-differentials
	\begin{align*}
		\sigma^2\alpha_1&\dotmapsto d\alpha_1\\
		\sigma^2 v_1&\dotmapsto dv_1.
	\end{align*}
	The class $\alpha_1$ is a permanent cycle since it must represent the image of $\alpha_1\in\pi_\ast j/(p,v_1)$ along the unit map, and the divided power classes $(d\alpha_1)^{(k)}$ are permanent cycles because they are in weight $0$, and there are no classes of weight $>1$. Therefore, by the Leibniz rule, the $E_{2p-1}$-page is isomorphic to
	\[
		\mathbb F_p[\mu_2]\otimes\Lambda[\lambda_2, a, \alpha_1]\otimes\Gamma[\gamma_p(d\alpha_1)]
	\]
	where $\mu_2,\lambda_2$ and $a$ represent $(\mu_1)^p,(\sigma^2 v_1)^{p-1}dv_1$ and $(\sigma^2 \alpha_1)\gamma_{p-1}(d\alpha_1)$, respectively.
		
	For degree reasons, the only possible further nonzero differential is
	\[
		d_{p-1}(\lambda_2) \doteq \alpha_1
	\]
	To prove that this differential actually happens, it is enough to show that
	\[
	\pi_{2p^2-2}\THH(j)/(p,v_1)=0.
	\]
	By \Cref{cor:jconn}, there is a $(4p^2-4p-2)$-connective map $j \oplus \Sigma^2\fib(\SP_p \to j) \to \THH(j)$, so it suffices to show that
	\[
		\pi_{2p^2-2}(j/(p,v_1)) = \pi_{2p^2-2}(\Sigma^2\fib(1_j)/p,v_1)=0.
	\]
	The former group is clearly $0$. The latter is $0$ from the computation of the Adams--Novikov $E_2$-page for $\SP/(p,v_1)$ in low degrees (see the discussion after \cite[Theorem 4.4.9]{ravenelgreen} and Theorem 4.4.8 of op. cit.).

The last nontrivial differential of the spectral sequence is displayed for $p=3$ in \Cref{fig:thhj}.
	
	We now check for $p\geq 5$ that there are no multiplicative extension problems in our description of the commutative ring structure on $\pi_*\THH(j)/(p,v_1)$. If we choose $\gamma_{p^i}b$ to be detected by $(\gamma_{p^{i+1}}(d\alpha_1))$, the relations $\gamma_{p^i}(b)^p=0$ follow since there is nothing of higher filtration in that degree. Let $\mu_2$ be any lift of $(\sigma^2v_1)^p$. The homology of the CDGA $\Lambda_{\FF_p}[\alpha_1,\lambda_2,a], d(\lambda_2) = a\alpha_1$ is $6$-dimensional over $\FF_p$, given by $$\{1,a,\alpha_1,\lambda_2a,\lambda_2\alpha_1,\lambda_2a\alpha_1\}$$ Let $\alpha_1,x, y,z$ denote lifts of the classes $\alpha_1,a,\lambda_2a,\lambda_2\alpha_1$ respectively (so that $\alpha_1y$ is a lift of $\lambda_2a\alpha_1$). The relation $\alpha_1y=-xz$ holds because it is true on the associated graded and there is nothing of higher filtration in that degree. The classes $\alpha_1z,yz,x\alpha_1$ are $0$ because there are no nonzero classes in degree $(p+1)(2p-2),2p^2-1+2(2p-3),2(2p^2-1)+(2p-3)+ p(2p-2)+1$ respectively. The only remaining relation, $xy=0$, occurs because it happens on the associated graded, and there is nothing of higher filtration.
\end{proof}

\begin{rmk}
	For $p=3$, it is more complicated to figure out the multiplicative extensions, since the homotopy ring is not necessarily associative. Many of the multiplicative extensions can be ruled out using the Postnikov filtration on $j/(3,v_1)$, but not all of them: for example this doesn't rule out the possible non-associative extension $x(x \mu_2^2) = zb^2$ in degree $62$.
\end{rmk}

\section{THH of finite extensions}\label{sec:finiteextn}
	In this section, we shall make the analogous computations for the THH of $j_{\zeta,k}:=\ell_p^{hp^k\mathbb Z}$, $ju_{\zeta,k}$, and and also of $j_k:=\tau_{\geq0}j_{\zeta,k}$ for $p>2$, which are introduced as filtered rings in \Cref{defn:lkofilt}. $j_{\zeta,k}$ is a $\ZZ/p^k$ Galois extension of $j_{\zeta}$ in $\Sp_p$. The computations are very similar to the cases of $j_{\zeta}$ and $j$, so we shall only point out the differences from the proofs of those cases.

	\begin{thm}\label{thm:jzetafinextn}
		There is an isomorphism of rings for $p>2$
		\[
			\pi_\ast \THH(j_{\zeta,k})/(p,v_1) \simeq \pi_*(\THH(\ell_p)/(p,v_1))\otimes \Lambda[\zeta]\otimes \ctf{\ZZ_p}
		\]
		and for $p=2$
		\[
		\pi_\ast \THH(j_{\zeta,k})/(2,\eta,v_1) \simeq \pi_*(\THH(\ko_2)/(2,\eta,v_1))\otimes \Lambda[\zeta]\otimes \ctf{\ZZ_2}
		\]
		
		The maps $\THH(j_{\zeta,k})/(p,v_1) \to \THH(j_{\zeta,k+1})/(p,v_1)$ 
		on $\pi_*$ are the identity on the $\THH(\ell_p)/(p,v_1)$ component,
		send $\zeta$ to $0$, and are the restriction map $\ctf{\ZZ_p} \to \ctf{p\ZZ_p}\cong \ctf{\ZZ_p}$.
	\end{thm}
	\begin{proof}
		The proof strategy is the same as in \Cref{thm:jzetathh} and \Cref{thm:thhjzeta2}. One difference is that for $k\geq1$, the class $\lambda_1$ in the spectral sequence $\THH(j_{\zeta,k}^{\gr})/(p,v_1) \implies \THH(j_{\zeta,k})/(p,v_1)$ is a permanent cycle, which can be seen from the Leibniz rule. In particular, \Cref{rmk:lambda1} doesn't apply for $k\geq1$, and as noted in the remark, this difference doesn't affect the final answer.
		
		The claim about the maps $\pi_*\THH(j_{\zeta,k})/(p,v_1) \to \THH(j_{\zeta,k+1})/(p,v_1)$ can be deduced at the level of associated graded of the filtrations. For example, by choosing elements $\lambda_1,\lambda_2,\sigma^2v_2$ in $\THH(j_{\zeta})/(p,v_1)$, one sees that their images in $\THH(j_{\zeta,k})/(p,v_1)$ are valid generators of the corresponding classes. To see what the transition maps do on $\Lambda[\zeta]\otimes \ctf{\ZZ_p}$, we can use \Cref{lem:spheretriv} since these classes are in the image of $\THH(\SP_p^{B\ZZ})$. It then follows that map sends $\ctf{\ZZ_p} \to \ctf{p\ZZ_p}$ given by restriction of functions, and $\zeta$ goes to $p \zeta=0$ because that is what happens on the level of mod $p$ cohomology of the $p$-fold cover map $S^1 \to S^1$.
	\end{proof}
	
		We next explain the computation for $ju_{\zeta,k}$, which is nearly identical to that of $j_{\zeta,k}$
	\begin{thm}\label{thm:juk}
		For each $k\geq0$, there is an isomorphism of rings
		\[
		\pi_\ast \THH(ju_{\zeta,k})/(2,v_1) \simeq \pi_*(\THH(\ell_2)/(2,v_1))\otimes \Lambda[\zeta]\otimes \ctf{\ZZ_2}
		\]
		
		The maps $\THH(ju_{\zeta,k})/(p,v_1) \to \THH(ju_{\zeta,k+1})/(p,v_1)$ on $\pi_*$ are the identity on the
		$\THH(\ell_2)/(2,v_1)$ component, send $\zeta$ to $0$, and are the restriction map $\ctf{\ZZ_2} \to \ctf{2\ZZ_2}\cong \ctf{\ZZ_2}$
	\end{thm}
	
	\begin{proof}
		The proof is nearly exactly as the proof of \Cref{thm:jzetafinextn} for $p>2$. The only difference is that in checking multiplicative extension problems in spectral sequences, one must check that odd degree classes square to zero (since we are at the prime $2$). This always follows because the square lands in a zero group; see \Cref{fig:thhju} for a chart.
	\end{proof}

	Our argument to compute $\THH(j_k)$ for $k\geq1$ uses Dyer--Lashof operations to produce permanent cycles, so we first give $j_k/(p,v_1)$ an $\EE_{\infty}$-structure.
	
	\begin{prop}\label{prop:jkmodstuffeinfty}
		For $k\geq1$, $j_k/(p,v_1)$ admits the structure of an $\EE_{\infty}$-algebra under $j_k$ that is a trivial square zero extension of $\FF_p$ by $\Sigma^{2p-2}\FF_p$.
	\end{prop}
	
	\begin{proof}
		To construct $j_k/(p,v_1)$ as an $\EE_{\infty}$-ring, we first begin with $\tau_{\leq 2p-3}j_{k}$, whose homotopy groups are $\ZZ_p$ in degree $0$ and $\ZZ/p^{k+1}$ in degree $2p-3$, where $\alpha_1$ is a $p$-torsion class in degree $2p-3$.
		
		By \cite[Corollary 7.4.1.28]{HA} this is a square zero extension of $\ZZ_p$ by $\Sigma^{2p-3}\ZZ/p^{k+1}$, i.e it fits into a pullback square
		
		\begin{center}
			\begin{tikzcd}
				\tau_{\leq2p-3}j_k \ar[r]\ar[d] &\ZZ_p \ar[d]\\
			\ZZ_p\ar[r] & \ZZ_p\oplus \Sigma^{2p-2}\ZZ/p^{k+1}
			\end{tikzcd}
		\end{center}
	
		By using the map $\ZZ/p^{k+1} \to \ZZ/p$ that kills every multiple of $p$ (including $\alpha_1$ since $k\geq1$), we can produce an $\EE_{\infty}$-algebra $R$ under $\tau_{\leq 2p-3}j_k$ defined as the pullback
		
		\begin{center}
			\begin{tikzcd}
				R \ar[r]\ar[d] &\ZZ_p \ar[d]\\
				\ZZ_p\ar[r] & \ZZ_p\oplus \Sigma^{2p-2}\ZZ/p
			\end{tikzcd}
		\end{center}
	
		We claim that $R$ is a trivial square zero extension of $\ZZ_p$. To see this, square zero extensions of $\ZZ_p$ by $\Sigma^{2p-1}\FF_p$ are classified by maps of $\ZZ_p$-modules $L_{\ZZ_p/\SP_p} \to \Sigma^{2p-1}\FF_p$, where $L_{\ZZ_p/\SP_p}$ denotes the $\EE_{\infty}$ relative cotangent complex.  By \cite[Theorem 7.4.3.1]{HA}, since $\SP_p \to \ZZ_p$ is $2p-3$-connective, there is a $4p-4$-connective map 
		
		$$\ZZ_p \otimes_{\SP_p}\cof(\SP_p \to \ZZ_p) \to L_{\SP_p/\ZZ_p}$$
		
		showing that $\pi_{2p-2}L_{\ZZ_p/\SP_p}$ is $\FF_p$. It follows that up to isomorphism, there is a unique nontrivial square zero extension of $\ZZ_p$ by $\Sigma^{2p-3}\FF_p$. But $\tau_{\leq2p-3}\SP_p$ must be this nontrivial extension, since $\alpha_1\neq0$ there. Since $\alpha_1=0$ in $R$, it follows that $R$ is the trivial square zero extension $\ZZ_p \oplus \Sigma^{2p-3}\FF_p$. Thus $\tau_{\leq2p-3}(R\otimes_{\ZZ_p}\FF_p)$ is an $\EE_{\infty}$-$\FF_p$-algebra under it that is a trivial square zero extension of $\FF_p$ by $\Sigma^{2p-2}\FF_p$. Since $R$, and $v_1=0$ modulo $p$, there is a unital map $j_k/(p,v_1) \to R$, which is an equivalence because it is on homotopy groups.
	\end{proof}
		
	\begin{thm}
		For $k\geq1,p>2$, there is an isomorphism
		\[
			\pi_\ast\THH(j_k)/(p,v_1) \simeq \pi_*\THH(\ell_p)/(p.v_1) \otimes \Lambda[\alpha_{1/p^k}]\otimes\Gamma[d\alpha_{1/p^k}]
		\]
		where $|\alpha_{1/p^k}| = 2p-2$ and $|\sigma\alpha_{1/p^k}| = 2p-1$.
	\end{thm}
	\begin{proof}
		The proof of \Cref{prop:thhjgr} carries over exactly for $j_k$ to give an isomorphism
		$$	\pi_*\THH(j_k^{\gr})/(p,v_1) \cong \pi_*\THH(\ZZ_p)/p\otimes \pi_*\HoH(\tau_{\geq0}\FF_p[v_1]^{Bp^k\ZZ}/\FF_p)/v_1	$$
		The second tensor factor on the right hand side by \Cref{prop:hhjgr/p} is $\Lambda[\alpha_{1/p^k},dv_1]\otimes \Gamma[d \alpha_{1/p^k}]$\footnote{As an algebra this doesn't depend on $k$, but we have given names depending on $k$ to indicate that the exterior class $\alpha_{1/p^k}$ is sent to $0$ in $\THH(j_{k+1}^{\gr})/(p,v_1)$.}. 
		
		In the spectral sequence for $\THH(j_k^{\gr})/(p,v_1) \implies \THH(j_k)/(p,v_1)$, there is a differential $d_{2p-2}\sigma^2v_1 = d v_1$ arising as in \Cref{thm:thhj}, but the target of the differential from $\sigma^2\alpha_1$, which is $\sigma\alpha_1$, is zero since $\alpha_1 = 0$ in $j_k/(p,v_1)$. In fact, the class $\sigma^2\alpha_1$ is a permanent cycle since it can be constructed using a nullhomotopy of $\alpha_1$. Let $\lambda_1$ be a class in $\THH(j_k)/(p,v_1)$ detecting this.  
		
		By \Cref{prop:jkmodstuffeinfty}, $j_{k}/(p,v_1)$ is an $\EE_{\infty}$-algebra under $j_k$ that is an $\EE_{\infty}$-$\FF_p$-algebra, so $\THH(j_k)/(p,v_1) \cong \THH(j_k)\otimes_{j_k}j_k/(p,v_1)$ is an $\EE_{\infty}$-$\FF_p$-algebra with Dyer--Lashof operations. We define $\lambda_2$ to be the $\EE_2$-Dyer--Lashof operation on $\lambda_1$. In $\THH(\ell_p)/(p,v_1)$, this operation on the class $\lambda_1$ gives the class $\lambda_2$ in $\pi_{2p^2-1}\THH(\ell_p)/(p,v_1)$ \cite[Section 2]{ausonirognes}, which is detected by $\sigma^2v_1^{p-1}dv_1$ in the spectral sequence for $\THH(\ell_p^{\fil})/(p,v_1)$ by \Cref{exm:THHell}. Since maps of filtered objects can only increase filtrations in which elements are detected, it follows that $\lambda_2$ must also be detected by $\sigma^2v_1^{p-1}dv_1$ in $\THH(j_k)/(p,v_1)$, so that class is a permanent cycle. The class $\alpha_{1/p^k}$ is a permanent cycle since it is in the image of the unit map, and the classes in $\Gamma[d\alpha_{1/p^k}]$ must be permanent cycles for degree reasons, so there are no further differentials. There are no even degree classes of positive weight, so classes representing the divided powers of $d\alpha_{1/p^k}$ have zero $p^{th}$-power for degree reasons. For degree reasons there can be no further multiplicative extensions.
	\end{proof}

\section{$\TC$ in the stable range}\label{sec:stable}

$\TC$ is an important invariant of rings, partially because of the Dundas--Goodwillie-McCarthy theorem \cite{DGM}, which says that for nilpotent extensions of rings, the relative $K$-theory is the relative $\TC$.

\begin{thm}[Dundas--Goodwillie--McCarthy]\label{thm:dgm}
	Let $f:R\to S$ an $i$-connective map of connective $\EE_1$-rings, for $i\geq1$. Then there is a pullback square
	\begin{center}
		\begin{tikzcd}
			K(R) \ar[r]\ar[d] & K(S) \ar[d]\\
			\TC(R) \ar[r] & \TC(R)
		\end{tikzcd}
	\end{center}
\end{thm}

A precursor to this theorem is a result of Waldhausen\footnote{Although Waldhausen proves this result for $\EE_1$-$\ZZ$-algebras, the proof works equally well for any $\EE_1$-algebra: see for example \cite[Proposition 3.3]{levy2022algebraic}.}, which computes the first nonvanishing homotopy group of $\fib\TC(f) \cong \fib K(f)$ in terms of Hochschild homology.

\begin{prop}[{Waldhausen \cite[Proposition 1.2]{waldhausen1978algebraic}}]\label{prop:tchurewicz}
	Let $f:R \to S$ be an $i$-connective map of connective $\EE_1$-algebras for $i\geq1$. Then $\fib(K(f)) \cong \fib(\TC(f))$ is $(i+1)$-connective, with $\pi_{i+1}\fib(K(f)) \cong \mathrm{HH}_0(\pi_0S;\pi_i\fib f)$.
\end{prop}

Our goal in this section is to refine \Cref{prop:tchurewicz} to compute the spectrum $\fib(K(f))$ in the stable range in terms of $\THH$. We use this to understand the maps $K(\SP_p) \to K(\ZZ_p)$ and $K(j_{\zeta}) \to K(\ZZ_p^{B\ZZ})$ in the stable range. 

Given a map of $\EE_1$-rings, $R \to S$, the relative $\EE_1$-cotangent complex $L_{S/R}$ is the $S$-bimodule given by the fiber of the multiplication map $S\otimes_RS \to S$\footnote{See for example \cite[Remark 7.4.1.12]{HA}.}. Our result is as follows:

\begin{thm}\label{thm:tcstablerange}
	Given a map of ring spectra $f:R \to S$, there is a natural map $\fib\TC(f) \to \THH(S;L_{S/R})$. If $f$ is an $n$-connective map of $-1$-connective rings for $n\geq 1$, this natural map is  $2n+1$-connective.
\end{thm}

\begin{rmk}
	In fact the map of \Cref{thm:tcstablerange} is the linearization map in the sense of Goodwillie calculus, of the functor $f \mapsto \fib (\TC(f))$. See \cite{hesselholt1994stable,dundas1994stable} for a variant of this, where one considers only trivial square-zero extensions of $S$ rather than arbitrary $\EE_1$-ring maps.
\end{rmk}

We first construct the natural transformation using the following lemma.

\begin{lem}\label{lem:cycinv}
	Let $f:R \to S$ be a map of $\EE_1$-rings. Then there is a natural equivalence $\THH(R;S) \cong \THH(S;S\otimes_RS)$ making the diagram below commute.
	
		\begin{center}
		\begin{tikzcd}
			\THH(R;S)\ar[rr]\ar[dr] & &\THH(S;S) \\
			& \ar[ur] \THH(S;S\otimes_RS)& 
		\end{tikzcd}
	\end{center}
\end{lem}

\begin{proof} 
	Consider the map $f^*:\Mod(R) \to \Mod(S)$ and its right adjoint $f_*:\Mod(S) \to \Mod(R)$. The composite $f^*f_*$ corresponds to the $S$-bimodule $S\otimes_RS$, and the composite $f_*f^*$ corresponds to the $R$-bimodule $S$. Since $\THH$ of a bimodule is the trace of the bimodule as an endomorphism in presentable stable categories, cyclic invariance of the trace gives the desired equivalence $\THH(R;S) \cong \THH(S;S\otimes_RS)$. There is a diagram
	
\[\begin{tikzcd}
	{\Mod(R)} && {\Mod(S)} \\
	{\Mod(S)} && {\Mod(S)}
	\arrow["{1_S}"', shift right=3, from=1-3, to=2-3]
	\arrow["{1_S}"', shift right=3,  from=2-3, to=1-3]
	\arrow["{f^*}"', shift right=3, from=1-1, to=2-1]
	\arrow["{f_*}"', shift right=3, from=2-1, to=1-1]
	\arrow["{f^*}"{description}, from=1-1, to=1-3]
	\arrow["{1_S}"{description}, from=2-1, to=2-3]
\end{tikzcd}\]
	where we use the natural transformation $\epsilon:f^*f_* \to 1_S$ and $1_{f^*}$ to fill in the $2$-morphisms in the diagram. We recall that by Morita theory, colimit preserving functors between module categories are uniquely given by tensoring with a bimodule. Thus the horizontal maps in the diagram induce at the level of bimodules the maps $f^*f_* \implies 1_S$ and $f_*f^* \implies 1_S$ which induce the maps $\THH(R;S),\THH(S;S\otimes_RS) \to \THH(S;S)$ in the triangle of the lemma statement. The $C_2$-action on $\THH(S)$ coming from writing $1_S$ as $1_S\circ 1_S$ corresponds to restricting the $S^1$-action on $\THH(S)$ to $C_2 \subset S^1$. It follows that the claimed diagram naturally commutes because $S^1$ is connected, so the rotation by $\pi$ action on $\THH(S)$ is homotopic to the identity.
\end{proof}

\begin{cnstr}\label{cnstr:nattrans}
	We now construct a map $\fib(\TC(f)) \to \THH(S;L_{S/R})$, which is a natural transformation when viewed as a map between functors $\Alg(\Sp)^{\Delta^1} \to \Sp$.
	
	The map $f:R \to S$ gives a natural map as follows: composing the map $\TC(R) \to \THH(R)$ with $\THH(R) \to \THH(R;S)$, we obtain a commutative square
	
	\begin{center}
		\begin{tikzcd}
			\TC(R)\ar[r]\ar[d] & \ar[d]\TC(S)\\
			\THH(R;S)\ar[r] & \THH(S;S)
		\end{tikzcd}
	\end{center}
	Taking horizontal fibers and using the isomorphism of \Cref{lem:cycinv}, we obtain the desired natural transformation.
\end{cnstr}

We will first prove \Cref{thm:tcstablerange} in the case $R \to S$ is a square-zero extension with ideal $M$. To do this, we consider the square-zero extension as a filtered $\EE_1$-ring with underlying $R$ and associated graded $S\oplus M[1]$. Then $\THH(R)$ is a filtered $S^1$-equivariant spectrum, and the Frobenius maps $\Phi_p:\THH(R) \to \THH(R)^{tC_p}$ send filtration $i$ to filtration $ip$, so in particular can be thought of as filtration preserving maps, since the filtration is only in nonnegative degrees.

The key input we use is the computation of $\THH$ of a trivial square-zero extension as an $S^1$-equivariant spectrum:
\begin{prop}[{\cite[Proposition 4.5.1]{raskin2018dundas}}]\label{prop:thhsq0}
	For $S\oplus M$ the trivial square-zero extension of an $\EE_1$-ring $S$ by a bimodule $M$, there is an $S^1$-equivariant graded equivalence $\THH(S\oplus M) \cong \THH(S) \oplus \bigoplus_{m=1}^{\infty}\Ind_{\ZZ/m\ZZ}^{S^1}\THH(S;(\Sigma M)^{\otimes m})$
\end{prop}

Here $\Ind_{\ZZ/m\ZZ}^{S^1}$ is the right adjoint of the forgetful functor from $S^1$-equivariant spectra to $\ZZ/m\ZZ$-spectra, and the $\ZZ/m\ZZ$-action on $\THH(S;(\Sigma M)^{\otimes m})$ comes from cyclically permuting the tensor factors.

We also record a key property of the $\THH$ of $-1$-connective rings that we use:

\begin{lem}\label{lem:thhconn}
	Let $R \to S$ be an $n$-connective map of $-1$-connective rings, and $M$ a connective $S$-bimodule. Then $\THH(S;M)$ is connective, and the map $\THH(R;M) \to \THH(S;M)$ is $n+1$-connective.
\end{lem}

\begin{proof}
	Both of these follow from examining the associated graded coming from the cyclic bar complex computing $\THH(R;M)$ and $\THH(S;M)$. For the latter is given by $\Sigma^mS^{\otimes m}\otimes M$ which indeed is connective, and $\Sigma^mS^{\otimes m}\otimes M \to \Sigma^mR^{\otimes m}\otimes M$ is $n+m$-connective for $m\geq1$ and an isomorphism for $m=0$.
\end{proof}

\begin{prop}\label{prop:sq0stable}
	Let $f:R\to S$ be an $n$-connective square-zero extension of $-1$-connective $\EE_1$-rings for $n\geq0$. Then the map $\fib\TC(f) \to \THH(S;L_{S/R})$ is $2n+1$-connective.
\end{prop}

\begin{proof}
	We consider the map $\fib\TC(f) \to \fib\THH(f) \to \THH(S;L_{S/R})$ as a map of filtered spectra, viewing $S$ as a filtered $\EE_1$-ring with associated graded $R\oplus M$. By \Cref{prop:thhsq0}, $\gr(\fib\THH(R)) \cong \bigoplus_{m=1}^{\infty}\Ind_{\ZZ/m\ZZ}^{S^1}\THH(S;(\Sigma M)^{\otimes m})$ as an $S^1$-spectrum. Since the Frobenius map is zero on associated graded since it takes filtration $i$ to $ip$, so we learn that $\gr_m(\fib \TC(f)) \cong (\Sigma \Ind_{\ZZ/m\ZZ}^{S^1}\THH(S;(\Sigma M)^{\otimes m}))_{hS^1}$\footnote{See also \cite[Theorem 4.10.1]{raskin2018dundas}.}. In particular, since $S$ is $-1$-connective and $n\geq 0$, the connectivity of these terms goes to $\infty$ as $m \to \infty$ via \Cref{lem:thhconn} so the filtration on $\TC$ is complete. Since $\Ind_{\ZZ/m\ZZ}^{S^1}$ decreases connectivity by $1$, we learn that $\gr_m(\fib \TC(f))$ is $(n+1)m-1$-connective. In particular, the map $\fib \TC(f) \to \gr_1\fib\TC(f)$ is $2n+1$-connective.
	
	To finish, it suffices to show the following two claims:
	\begin{enumerate}
		\item $\THH(S;L_{S/R}) \to \gr_1\THH(S;L_{S/R})$ is $2n+2$-connective.
		\item $\gr_1\fib\TC(f) \to\gr_1 \fib\THH(S;L_{S/R})$ is an isomorphism.
	\end{enumerate}
	The claim $(1)$ follows from the fact that $\gr L_{S/R} \cong L_{S/S\oplus M} \cong \bigoplus_{m=1}^{\infty}(\Sigma M)^{\otimes_S m}$, and $(\Sigma M)^{\otimes_S m}$ is $2n+2$-connective for $m\geq 2$.
	
	For claim $(2)$, we see that $$\gr_1\fib \TC(f) \cong \Sigma(\Ind_{\ZZ/1\ZZ}^{S^1}\THH(S;\Sigma M))_{hS^1}\cong (\Ind_{\ZZ/1\ZZ}^{S^1}\THH(S;\Sigma M))^{hS^1} \cong \THH(S;\Sigma M)$$ $\Sigma M$ is exactly $\gr_1L_{S/R}$, and $\THH(S;\gr_1L_{S/R}) \cong \gr_1\THH(S;L_{S/R})$ since $S$ is entirely in grading $0$, so we are done.
\end{proof}

We prove \Cref{thm:tcstablerange} by reducing to the case of a square-zero extension. First, we produce a natural way to factor a map of $\EE_1$-rings through a square-zero extension. We recall that given a $S'$-$S$-bimodule $M$ with a unit map $1:\SP \to M$, the pullback $S'\times_MS$ admits an $\EE_1$-algebra structure where the maps $S' \to M$ and $S \to M$ are the left $S'$-module and right $S$-module maps adjoint to the map $1$. This ring structure can be constructed as the endomorphism ring of the triple $(S',S,S \to S'\otimes_{S'}M)$ viewed as an object of the oplax limit $\Mod(S)\vec{\times}_{M}\Mod(S')$ (see \cite[Construction 2.5]{land2023k} and \cite[Construction 4.1]{burklund2021k}). When $M$ comes from a cospan of ring maps $S' \to R \leftarrow S$, this agrees with the pullback of the span of rings by \cite[Lemma 1.7]{land2019k}.

\begin{cnstr}
	Given a map $f:R \to S$, we consider $S\otimes_RS$ as an $S$-$S$-bimodule with unit $1$. We define $R_{f,2}$ to be the $\EE_1$-ring given by $S\times_{S\otimes_RS}S$.
\end{cnstr}

\begin{lem}\label{lem:univsq0}
	We have natural maps $R \xrightarrow{h} R_{f,2} \xrightarrow{g} S$. If $R \to S$ is an $n$-connective map of connective rings for $n\geq 0$, then $h$ is $2n$-connective, $g$ is $n$-connective, and $g$ is a square-zero extension.
\end{lem}

\begin{proof}
	The fiber of $h:R \to R_{f,2}$ is the total fiber of the square 
	\begin{center}
		\begin{tikzcd}
			R \ar[r]\ar[d] &S \ar[d]\\
			S\ar[r] & S\otimes_{R}S
		\end{tikzcd}
	\end{center}
	which is $\fib f\otimes_R\fib f$, which is $2n$-connective. Since $f$ is $n$-connective, it follows that $g$ is too. It remains to show that $g$ is a square-zero extension, which will follow if we identify $S\otimes_RS$ as an $S$-bimodule with unit with the associated structure on $S \oplus L_{S/R}$ coming from the cospan of rings $S \to S\oplus L_{S/R} \leftarrow{S}$ corresponding to the universal derivation. But since $R$ maps into the pullback of this cospan (since it is the universal square-zero extension of $S$ under $R$) we have a square of ring maps
	
	\begin{center}
		\begin{tikzcd}
			R\ar[r]\ar[d] &S \ar[d]\\
			S\ar[r] & S\oplus L_{S/R}
		\end{tikzcd}
	\end{center}

	which defines an isomorphism of unital $S$-bimodules $S\otimes_RS \to S\oplus L_{S/R}$.
\end{proof}

\begin{proof}[Proof of \Cref{thm:tcstablerange}]
	We consider the maps $h,g,f$ as in \Cref{lem:univsq0}, giving us the diagram
	\begin{equation}
	\begin{tikzcd}
		\fib\TC(h)\ar[r]\ar[d] &\fib\TC(f)\ar[r]\ar[d] & \fib\TC(g)\ar[d] \\
		\THH(R_{f,2};L_{R_{f,2}/R})\ar[r] &\THH(S;L_{R/S}) \ar[r] &\THH(S;L_{R_{2,f}/S})
	\end{tikzcd}
	\end{equation}
To produce a nullhomotopy of the composite of the lower horizontal maps, we identify them with the vertical fibers of the following cofiber sequence using \Cref{lem:cycinv}:

\begin{center}
	\begin{tikzcd}
		\THH(R;R_{f,2})\ar[r]\ar[d] & \THH(R;S)\ar[r]\ar[d] & \THH(R_{f,2};S)\ar[d] \\
		\THH(R_{f,2})\ar[r] &\THH(S) \ar[r] &\THH(S)
	\end{tikzcd}
\end{center}

The map $\THH(R_{f,2}) \to \THH(S)$ lifts to $\THH(R_{f,2};S)$, and this lifting provides the desired nullhomotopy. Moreover, we see that the fiber of the map $\THH(R_{f,2};L_{R_{f,2}/R}) \to \THH(S;\fib L_{R/S}\to L_{R_{2,f}/S})$ is identified with the total fiber of the square

\begin{center}
	\begin{tikzcd}
		\THH(R;R_{f,2})\ar[r]\ar[d] & \THH(R;S)\ar[d]\\
		\ar[r]\THH(R_{f,2}) & \THH(R_{f,2};S)
	\end{tikzcd}
\end{center}
	which is the fiber of the map $\THH(R;\fib g) \to \THH(R_{f,2};\fib g)$. By \cite[Lemma 3.2]{levy2022algebraic}, since $h$ is $2n$-connective and $\fib g$ is $n$-connective, we see that this map is $3n+1$-connective.
	
	We next observe that in the right square of diagram $(4)$, we know all maps except possibly the vertical map which we want to show is $2n+1$-connective. Indeed, $\fib\TC(h)$ is $2n+1$-connective by \Cref{lem:univsq0} and \Cref{prop:tchurewicz}, the right vertical map is $2n+1$-connective by \Cref{prop:sq0stable}, and the lower horizontal map is $2n+1$-connective since the map $S\otimes_RS \to S\otimes_{R_{f,2}}S$ is $2n+1$-connective by \cite[Lemma 3.2]{levy2022algebraic}. It follows that the middle vertical map in diagram $(4)$ is $2n$-connective. But since $f$ is an arbitrary $n$-connective map and $h$ is $2n$-connective, we learn that the left vertical map is $4n$-connective. It follows that the middle vertical map is $2n+1$-connective since it is an extension of a $2n+1$-connective map and a $4n$-connective map since $n\geq1$.
\end{proof}

\begin{rmk}
	There is a version of \Cref{thm:tcstablerange} for a $0$-connective map of connective rings, but one must ask that $\pi_0R\to \pi_0S$ has a nilpotent kernel.
\end{rmk}

\subsection{Applications to the sphere and the $K(1)$-local sphere}

We now apply \Cref{thm:tcstablerange} to the map $\SP_p \to \ZZ_p$ for $p\geq 2$ to understand the map $\TC(\SP_p) \to \TC(\ZZ_p)$ in the stable range. The proposition below contains a key ingredient of \cite[Section 9]{bokstedt1993topological} used to understand the homotopy type of $\TC(\ZZ_p)$.

\begin{prop}\label{prop:spherestablerange}
	For $p>2$, the map $\pi_*\TC(\SP_p) \to \pi_*\TC(\ZZ_p)$ in degrees $\leq 4p-6$ is an isomorphism in all degrees except $2p-1$, where it is the map $p\ZZ_p \to \ZZ_p$.
\end{prop}

\begin{proof}
	By \Cref{thm:stablerange} we have a $4p-5$-connective map
	$$ \fib(\TC(\SP_p) \to \TC(\ZZ_p))\to \fib(\THH(\SP_p;\ZZ_p) \to \THH(\ZZ_p))$$
	
	The target of the map is $\fib(\ZZ_p \to \THH(\ZZ_p))$, which after applying $\tau_{\leq4p-5}$ is $\Sigma^{2p-2}\FF_p$. Thus it follows that there is a cofiber sequence $$\Sigma^{2p-2}\FF_p \to  \tau_{\leq4p-4}\TC(\SP_p) \to \tau_{\leq4p-4}\TC(\ZZ_p)$$ Recall that $\TC(\SP_p) \cong \SP_p \oplus \Sigma (\CC\PP^{\infty}_{-1})_p$ \cite{bokstedt1993cyclotomic}\footnote{see also \cite[Theorem 8.4]{krause2018lectures}}, and that $\pi_*\TC(\ZZ_p)/(p,v_1)$ is $\FF_p$ in odd degrees between $-1$ and $2p-1$, and in degrees $0,2p-2$, and $0$ in all other degrees \cite{bokstedt1993topological}\footnote{This argument is not circular, because $\TC(\ZZ_p)/(p,v_1)$ is computed without knowing this proposition.}. From this description, it follows that both $\TC(\SP_p)/(p,v_1)$ and $\TC(\ZZ_p)/(p,v_1)$ are $\FF_p$ in degrees $2p-2,2p-1$. Thus in the cofiber sequence above mod $(p,v_1)$, the class in degree $2p-2$ must go to $0$ and the class in degree $2p-1$ must go to the generator. It follows that integrally, the class must go to $0$, and that it maps to the $\ZZ_p$ in $\TC_{2p-1}(\SP_p)$ via the $p$-Bockstein, giving the conclusion.
\end{proof}

\begin{cor}\label{cor:jstablerange}
	For $p>2$, \Cref{prop:spherestablerange} also holds for $j$. In particular, the obstruction to lifting $\lambda_1 \in \TC(\ZZ_p)$ to $\TC(j)$ is up to a unit in $\FF_p$ the class $\sigma\alpha_1$ in $\THH(\ZZ_p;L_{\ZZ_p/j})$.
\end{cor}
\begin{proof}
	Since the map $\SP_p \to j$ is $2p^2-2p-2$-connective (see \Cref{cor:jconn}), the map $\SP_p \to \ZZ_p$ agrees with the map $j_p\to \ZZ_p$ in the stable range, so the analysis in \Cref{prop:spherestablerange} applies for $j$. In particular, the obstruction to lifting the class $\lambda_1 \in \TC(\ZZ_p)$ to $j$ is nonzero in $\THH(\ZZ_p;L_{\ZZ_p/j})$, so must be $\sigma\alpha_1$ up to a unit in $\FF_p$, since $\pi_{2p-2}\THH(\ZZ_p;L_{\ZZ_p/j}) \cong \FF_p$ is generated by this class.
\end{proof}

We now apply \Cref{thm:tcstablerange} to the map $j_{\zeta} \to \ZZ_{\zeta}$, and then make deductions about $K(L_{K(1)}\SP)$ in the stable range.

\begin{lem}\label{lem:jzetarelcot}
	There is an isomorphism $\Sigma^{2p-2}\FF_p \cong L_{\ZZ_{\zeta}/j_{\zeta}}$, where the generator is $\sigma(\alpha_1)$.
\end{lem}

\begin{proof}
	In fact, we claim that $L_{\ZZ_{\zeta}^{\gr}/j_{\zeta}^{\gr}} \cong \Sigma^{2p-2,0}\FF_p$ on the class $\sigma(\alpha_1)$ which implies the result, since this is the associated graded of $L_{\ZZ_{\zeta}/j_{\zeta}}$. To see this, we note that $L_{\ZZ_{\zeta}^{\gr}/j_{\zeta}^{\gr}}/p \cong L_{\ZZ_{\zeta}^{\gr}/p/j_{\zeta}^{\gr}/p}$. Since $j_{\zeta}^{\gr}/p \to \ZZ_{\zeta}/p$ is the augmentation of a polynomial algebra over the target on the class $v_1$, $L_{\ZZ_{\zeta}^{\gr}/p/j_{\zeta}^{\gr}/p} \cong \Sigma^{2p-1} \ZZ_{\zeta}^{\gr}/p$, where the generating class is $\sigma(v_1)$. In $j_{\zeta}^{\gr}$, there is a $p$-Bockstein differential $d_1v_1 = v_1\zeta = \alpha_1$, so applying the map $\sigma$, we get that $\sigma(v_1)$ has a $p$-Bockstein $d_1$-differential hitting $\zeta \sigma(v_1) = \sigma(\alpha_1)$. Thus we can conclude.
\end{proof}

The following proposition gives a way in which $\TC(j_{\zeta})$ does not behave as if the action on $\ell_p$ is trivial.
\begin{prop}\label{prop:jzetalambda1}
	For $p>2$, the image of the class $\lambda_1 \in \TC(\ZZ_p)/(p,v_1)$ in $\TC(\ZZ_{\zeta})/(p,v_1)$ does not lift to $\TC(j_{\zeta})/(p,v_1)$. The same statement is true for $K$-theory replacing $\TC$.
	
	\begin{proof}
		The result for $K$-theory is equivalent to the one for $\TC$ by \cite{levy2022algebraic}.
		We have a commutative square of maps
		
		\begin{center}
			\begin{tikzcd}
				\fib(\TC(j) \to \TC(\ZZ_p) \ar[r]\ar[d] &\fib(\TC(j_{\zeta})\to \TC(\ZZ_p^{B\ZZ}))\ar[d]\\
				\THH(\ZZ_p;L_{\ZZ_p/j})\ar[r] & \THH(\ZZ_{\zeta};L_{\ZZ_{\zeta}/j_{\zeta}})
			\end{tikzcd}
		\end{center}
		
		where the vertical maps are $4p-5$-connective by \Cref{thm:tcstablerange}.
		
		The lower horizontal map sends $\sigma(\alpha_1)$ to $\sigma(\alpha_1)$, the generator of $\pi_{2p-2}\THH(\ZZ_{\zeta};L_{\ZZ_{\zeta}/j_{\zeta}})$. But $\sigma(\alpha_1)$ since the class is the obstruction to lifting $\lambda_1$ from $\TC(\ZZ_p)$ to $\TC(\SP_p)$, we learn that the obstruction to lifting $\lambda_1$ from $\TC(\ZZ_p^{B\ZZ})$ to $\TC(j_{\zeta})$ is nontrivial. We also see that this obstruction is nonzero modulo $(p,v_1)$.
	\end{proof}
\end{prop}

\begin{thm}
	For $p>2$, there are isomorphisms $$\tau_{\leq 4p-6}\fib(\TC(j_{\zeta}) \to \TC(\ZZ_{\zeta})) \cong \Sigma^{2p-2}\ctf{\ZZ_p}$$ and $$K_{*}L_{K(1)}\SP \cong K_{*-1}\FF_p \oplus K_{*}\SP_p \oplus \pi_*\Sigma^{2p-2}\ctf{\ZZ_p}/\FF_p, \;\;*\leq 4p-6$$
\end{thm}

\begin{proof}
	The map $f:j_{\zeta} \to \ZZ_{\zeta}$ is $2p-3$-connective, so we learn that $\fib\TC(f) \to \THH(\ZZ_{\zeta};L_{\ZZ_{\zeta}/j_{\zeta}})$ is $4p-5$-connective using \Cref{thm:tcstablerange}. For the first statement, it suffices to show that $\tau_{\leq 4p-4}\THH(\ZZ_{\zeta};L_{\ZZ_{\zeta}/j_{\zeta}}) \cong \Sigma^{2p-2}\ctf{\ZZ_p}$. But using \Cref{cor:trivfixed} and \Cref{lem:jzetarelcot}, we learn
	$$\THH(\ZZ_{\zeta};L_{\ZZ_{\zeta}/j_{\zeta}}) \cong \THH(\ZZ_{\zeta};\Sigma^{2p-2}\ZZ_{\zeta}/p\otimes_{\SP_p^{B\ZZ}}\SP_p)$$ $$\cong \Sigma^{2p-2}\THH(\ZZ_{\zeta})/p\otimes_{\SP_p^{B\ZZ}}\SP_p \cong \Sigma^{2p-2}\THH(\ZZ_p)/p\otimes_{\FF_p}\ctf{\ZZ_p}$$
	
	Since $\pi_*\THH(\ZZ_p)/p$ is by \Cref{exm:THHZ} $\FF_p[\sigma^2\alpha_1,\sigma^2v_1]$, we indeed learn the claim.
	
	To get the statement about $K$-theory, by \cite{levy2022algebraic}, we have $K_*(L_{K(1)}\SP) \cong K_*(j_{\zeta}) \oplus K_{*-1}(\FF_p)$ and a map of cofiber sequences
		\begin{center}
	\begin{tikzcd}
		\fib(\TC(\SP_p) \to \TC(\ZZ_p))\ar[d] \ar[r]& K(\SP_p)\ar[d] \ar[r]& K(\ZZ_p)\ar[d]\\
		\fib(\TC(j_{\zeta}) \to \TC(\ZZ_{\zeta})) \ar[r]& K(j_{\zeta}) \ar[r]& K(\ZZ_p)
	\end{tikzcd}
\end{center}
where the fiber terms are $0$ after inverting $p$. The third terms after $p$-completion are $K_{2p-1}(\ZZ_p)_p \cong \ZZ_p$, generated by $\lambda_1$.
	
	We can split $\pi_*\Sigma^{2p-2}\ctf{\ZZ_p}$ into $\pi_*\Sigma^{2p-2}\ctf{\ZZ_p}/\FF_p \oplus \FF_p$ via the augmentation to $\FF_p$ coming from evaluation at $0$. Note that the image of the boundary map of the cofiber sequence on mod $p$ homotopy groups is the $\FF_p$-summand by \Cref{prop:jzetalambda1}. It follows that after quotienting out by $\Sigma^{2p-2}\ctf{\ZZ_p}/\FF_p$ in the homotopy groups of the first two terms in the lower cofiber sequence, it agrees with that of the upper cofiber sequence in degrees $\leq 4p-6$.

It follows that there is a short exact sequence for $*\leq 4p-6$:	$$0 \to \pi_*\Sigma^{2p-2}\ctf{\ZZ_p}/\FF_p \to K_*(j_{\zeta}) \to K_*(\SP_p) \to 0$$
	But the map $K(\SP_p) \to K(j_{\zeta})$ clearly splits this sequence, giving the result.
\end{proof}

\section{The Segal conjecture}\label{sec:seg}
	
	The Segal conjecture for a cyclotomic spectrum $X$ is the statement that the cyclotomic Frobenius map $X \to X^{tC_p}$ is an isomorphism in large degrees. Knowing the Segal conjecture for $\THH(R)\otimes V$ where $V$ is a finite spectrum is a key step in proving the Lichtenbaum--Quillen conjecture for $X$, i.e the fact that $\TR(X)$ (and hence $\TC(X)$) is bounded (see \cite{hahn2022redshift}). 
	
	Asking that the Segal conjecture hold for $\THH(R)\otimes V$  is a regularity and finiteness condition on $R$: for example it holds when $V$ is $p$-torsion and $R$ is a $p$-torsion free excellent regular noetherian ring with the Frobenius on $R/p$ a finite map \cite[Corollary 1.5]{mathew2021k}. In this section, we show that the Segal conjecture does hold for $j_{\zeta}$ for $p>2$ as well as the extensions $j_{\zeta,k}$, but doesn't hold for the connective covers $j$ and $j_k$. In particular the Lichtenbaum--Quillen conjecture doesn't hold for $j_k$, and our result is used in \cite{telescope} to show that it does hold for $j_{\zeta,k}$ for $p>2$.
	
	A related regularity phenomenon was noted in \cite{levy2022algebraic}, namely that $j_{\zeta}$ is regular\footnote{See \cite{ncg} for a discussion of regularity in the setting of prestable $\infty$-categories.} at the height $2$-locus: i.e the $t$-structure on $\Mod(j_{\zeta})$ restricts to a bounded $t$-structure on $\Mod(j_{\zeta})^{\omega}\otimes \Sp_{\geq2}$. This $t$-structure is the key point in relating $j_{\zeta}$'s algebraic $K$-theory to that of the $K(1)$-local sphere. On the other hand, $j$ is not regular at the height $2$-locus which is why its integral $K$-theory is not closely related to that of the $K(1)$-local sphere.

Our first goal is to show that for odd $p$, $j_{\zeta,k}$ satisfies the Segal conjecture. A key input is the proposition below, the proof of which is the same as in the reference, though the statement is somewhat more general.

We first recall as in \cite[Section C.5]{hahn2022redshift} that given a filtered $\ZZ^m$-graded $\EE_1$-ring $R^{\fil}$, the cyclotomic Frobenius map refines to a filtered map $$\varphi: L_p\THH(R^{\fil})\to \THH(R^{\fil})^{tC_p}$$ where $L_p$ is the operation on filtered spectra scaling the filtration and the gradings on $R$ by $p$.

\begin{prop}{\cite[Proposition 4.2.2]{hahn2022redshift}}\label{prop:segpoly}
	Let $R$ be an $\EE_1$-ring, and consider the $\ZZ^m$-graded polynomial algebra $R[a_1,\dots,a_n]:=R\otimes\bigotimes_1^n \SP[a_i]$, where each $a_i$ has positive weight\footnote{i.e it is nonnegative weight in each copy of $\ZZ$ in $\ZZ^m$, and positive weight in some copy of $\ZZ$.} and is even topological degree and $\SP[a_i]$ is the free $\EE_1$-algebra. The map
 \[
 \varphi: L_p\THH(R[a_1,\dots,a_n]) \to \THH(R[a_1,\dots,a_n])^{tC_p} 
 \]
 at the level of $\pi_*$
 is equivalent to the map
	$$\pi_*\THH(R)[a_i]\otimes \Lambda[da_i] \to \pi_*\THH(R)^{tC_p}[a_i]\otimes \Lambda[da_i]$$
	where the $a_i,da_i$ are sent to themselves. If $R$ is an $\EE_2$-algebra and $\SP[a_i]$ are given the $\EE_2$-algebra structures coming from \cite{rotation}, this is a homomorphism of rings.
\end{prop}

The following lemma is used to reduce showing the Segal conjecture is true to the associated graded of a filtration on the ring.

\begin{lem}\label{lem:segassgr}
	Let $C$ be a presentably symmetric monoidal stable category with a complete $t$-structure compatible with filtered colimits, and suppose that $f:R^{\fil} \to R'^{\fil}$ is a map of homotopy associative filtered rings in $C$, where the filtration on the source and target is complete.
	
	If there is an element $x \in \pi_*R := \pi_*\map(\unit,R), *>0$ such that the associated graded map $R^{\gr} \to R'^{\gr}$ is $n$-coconnective in the constant $t$-structure and sends a class detecting $x$ to a unit, then the map $R \to R'$ is also $n$-coconnective, and is equivalent to the map
	
	$$R \to R[x^{-1}]$$
\end{lem}

\begin{proof}
	First, since the filtrations are complete and the map $f$ is $n$-coconnective on associated graded, we learn that the fiber is $n$-coconnective on associated graded, and complete, so the underlying object is $n$-coconnective. 
	
	Let $\tilde{x}$ be an element in $\pi_{**}R^{\fil}$ whose underlying element is $x$ that is sent to a unit in $R'^{gr}$. Since the filtration on $R'$ is complete, it follows that $\tilde{x}$ is sent to a unit, which allows us to build a map $R^{\fil}[\tilde{x}^{-1}]\to R'^{\fil}$ via the colimit of the diagram
	
\begin{center}
\begin{tikzcd}
	{\Sigma^{|x|} R^{\fil}}\ar[r] & {\Sigma^{|x|}R'^{\fil}} \\
	R^{\fil}\ar[r] & {R'^{\fil}} \\
	{{}} & {{}}
	\arrow["{...}"{marking}, shift left=1, draw=none, from=2-1, to=3-1]
	\arrow["{...}"{marking}, shift left=1, draw=none, from=2-2, to=3-2]
	\arrow["x", from=1-1, to=2-1]
	\arrow["x", from=1-2, to=2-2]
\end{tikzcd}
\end{center}
Note that the horizontal maps become more and more coconnective and the right vertical maps are all equivalences. Then because the $t$-structure is complete and compatible with filtered colimits, we learn that in the colimit the map is an equivalence. We also learn that the filtration on $R^{\fil}[x^{-1}]$ is complete, allowing us to conclude.
\end{proof}

\begin{thm}[Segal conjecture for $j_{\zeta,k}$]\label{thm:seg}
	For $p>2$ and $k\geq0$, the map $\THH(j_{\zeta,k})/(p,v_1) \to \THH(j_{\zeta,k})^{tC_p}/(p,v_1)$ has $2p-3$-coconnective fiber, and is equivalent to the map 
	$$\THH(j_{\zeta,k})/(p,v_1) \to \THH(j_{\zeta,k})[\mu^{-1}]/(p,v_1)$$ where $\mu \in \pi_{2p^2}\THH(j_{\zeta,k})$.
\end{thm}
\begin{proof}
	Using the filtration on $j_{\zeta,k}$ constructed in \Cref{sec:filtration}, we get a filtered map $$\varphi: L_p\THH(j_{\zeta,k})/(p,\tilde{v_1}) \to \THH(j_{\zeta,k})^{tC_p}/(p,\varphi\tilde{v_1})$$
	
	By the proof of \Cref{thm:jzetathh} and \Cref{thm:jzetafinextn}, the class $\mu$ is detected in the spectral sequence for $\THH(j_{\zeta,k})/(p,v_1)$ by $(\sigma^2v_1)^p$. Thus by applying \Cref{lem:segassgr} for $C = \Sp$ and $R^{\fil} \to R'^{\fil}$ the maps in question, it suffices to show
	\begin{enumerate}[label = (\alph*) ]
		\item The filtration on the source and target are complete.
		\item The associated graded map inverts the class $\sigma^2v_1$ and is $2p-3$-coconnective.
	\end{enumerate}
	
	To see $(a)$, the source is complete by \Cref{lem:postnikovcomplete}. The Tate construction $(-)^{tC_p}$ sits in a cofiber sequence up to shifts between the orbits $(-)_{hC_p}$ and fixed points $(-)^{hC_p}$, so it suffices to show each of those is complete. The orbits are complete for connectivity reasons: in any finite range of degrees, the orbits are computed via a finite colimit. The fixed points are complete because complete objects are closed under limits.
	
	We turn to proving $(b)$. We further filter $j_{\zeta,k}^{\gr}$ by the $p$-adic filtration as $j_{\zeta,k}^{\gr}\otimes \ZZ_p^{\fil}$ and consider the map of filtered graded $\EE_{\infty}$-rings $L_{p}\THH(j_{\zeta,k})/(\tilde{p},v_1) \to \THH(j_{\zeta,k})^{tC_p}/(\varphi\tilde{p},\varphi v_1)$. We claim:
	
	\begin{enumerate}[label = (\roman*) ]
		\item The filtration on the source and target are complete.
		\item The associated graded map inverts the class $\sigma^2p$ and is $2p-3$-coconnective.
	\end{enumerate}

	Given these claims, the proof is complete, since $\sigma^2v_1$ is detected in the spectral sequence by $(\sigma^2p)^p$ (see \Cref{exm:THHZ}), so claim (b) follows from \Cref{lem:segassgr}.

	(i) follows from an argument identical to the argument for (a), the only difference being that we use \Cref{lem:padiccomplete} to see that the filtration on $\THH(j_{\zeta}^{\gr}\otimes \ZZ_p^{\fil})/(\tilde{p},v_1)$ is complete. To see (ii), by \Cref{lem:modp} the associated graded algebra is $\FF_p[v_0,v_1]^{B\ZZ}$, where the action is trivial. By \Cref{lem:grdmodp} we have $\pi_*\THH(\FF_p[v_0,v_1]^{B\ZZ})/(v_0,v_1) \cong \ctf{\ZZ_p}\otimes \Lambda[dv_0,dv_1,\zeta]\otimes \FF_p[\sigma^2p]$, where $|dv_0| = 1, |\zeta| = -1, |dv_1|= 2p-1$. It follows that if the Frobenius map mod $(v_0,v_1)$ inverts $\sigma^2p$, it is $2p-3$-coconnective, since it is injective on $\pi_*$, and an element in the cokernel of largest degree is $(\sigma^2p)^{-1}\sigma v_1\sigma v_0$, which is in degree $2p-2$.
	
	Thus it remains to see that the Frobenius map mod $(v_0,v_1)$ on $\pi_*$ inverts the class $\sigma^2p$. Since $\THH$ is a localizing invariant and $\SP^{B\ZZ}$ is a trivial square-zero extension as an $\EE_1$-algebra, by \cite[Theorem 4.1]{land2023k} we have a pullback square of bigraded $\THH(\FF_p)$-modules in cyclotomic spectra
	
	\begin{center}
		\begin{tikzcd}
			{\THH(\FF_p[v_0,v_1]^{B\ZZ}}) \ar[r]\ar[d] &\THH(\FF_p[v_0,v_1]) \ar[d]\\
			\ar[r] \THH(\FF_p[v_0,v_1]) & \THH(\FF_p[v_0,v_1][x_0])
		\end{tikzcd}
	\end{center}

	where $x_0$ is a polynomial generator in degree $0$. It thus suffices to show that for $$\THH(\FF_p[v_0,v_1][x_0]),\THH(\FF_p[v_0,v_1])$$ the cyclotomic Frobenius map inverts $\sigma^2p$. These statements follow from \Cref{prop:segpoly} with $R = \FF_p,\FF_p[x_0]$, using the Segal conjecture for these discrete rings which is well known: for example \cite[Corollary 1.5]{mathew2021k} implies the Frobenius is an isomorphism in large degrees, but since it sends $\sigma^2p$ to a unit \cite[Corollary IV.4.13]{nikolaus2018topological}, it must just invert $\sigma^2p$.
\end{proof}

\begin{rmk}
	The bound $2p-3$ in \Cref{thm:seg} is optimal: the map is injective on $\pi_*$, and a class of largest degree not in the image is $\mu^{-1}\lambda_1\lambda_2$, in degree $2p-2$,
\end{rmk}

Now we show that the Segal conjecture fails for $\THH(j_k)$.
	
\begin{thm}\label{thm:segfailj}
	For $p>2$ and $k\geq0$, the fiber of the Frobenius map $\THH(j_k)/(p,v_1) \to \THH(j_k)^{tC_p}/(p,v_1)$ is not bounded above. Thus $j_k$ does not satisfy the Lichtenbaum--Quillen conjecture, i.e $\TR(j_k)\otimes V$ is not bounded above for $V$ a finite type $3$ spectrum.
\end{thm}
	
	\begin{proof}
		
		First we note that the failure of the Segal conjecture implies the failure of the Lichtenbaum--Quillen conjecture by \cite[Proposition 2.25]{antieau2021cartier}, so we show that the Segal conjecture fails.
		
		We first show that $\mu \in \THH(j_k)/(p,v_1)$ is sent to a unit in $\THH(j_k)^{tC_p}/(p,v_1)$. It follows from the spectral sequences used to calculate $\THH(j_k)/(p,v_1)$ that the image of $\mu$ in $\THH(\FF_p)$ is $(\sigma^2p)^{p^2}$ up to a unit, which is sent under the Frobenius map to a class detected up to a unit by $t^{-p^2}$ in the Tate spectral sequence for $\THH(\FF_p)^{tC_p}/(p,v_1)$ by \cite[Corollary IV.4.13]{nikolaus2018topological}. This is the lowest filtration of the Tate spectral sequence in $\pi_{2p^2}$, so since in that filtration the map $\THH(j_k)/(p,v_1) \to \THH(\FF_p)/(p,v_1)$ is the map $\ZZ_p \to \FF_p$, we learn that the image of $\mu$ must be detected by a unit multiple of $t^{-p^2}$ in the Tate spectral sequence for $\THH(j_k)^{tC_p}$ and hence be a unit.
		
		If the Frobenius map has an element $x$ in the kernel, then $x\mu^i$ is also in the kernel for each $i$, so the fiber isn't bounded above. On the other hand, if the Frobenius map is injective, then the classes $\varphi(\mu)^{-1}\varphi((\sigma \alpha_{1/p^k})^{(pi)})$ are an infinite family of classes of increasing degree in $\THH(j_k)^{tC_p}$ that are not in the image of $\varphi$, so in this case too, we learn that the fiber is not bounded above.
	\end{proof}

\begin{rmk}
	In fact, $\pi_*\THH(j_k)^{tC_p}/(p,v_1)$ under the Frobenius map is the completion of $\pi_*\THH(j_k)[\mu^{-1}]/(p,v_1)$ at the ideal generated by $(\sigma\alpha_{1/p^k}^{(pi)})$ for each $i$, and the map is in particular injective on $\pi_*$.
\end{rmk}

\bibliographystyle{alpha}
\bibliography{ref}
\end{document}